\documentclass[a4paper,11pt]{article}

\usepackage[utf8]{inputenc}

\usepackage[numbers,sort&compress]{natbib}
\usepackage{amsmath,amssymb,amsfonts}
\usepackage{listings}
\usepackage{amsbsy}
\usepackage{color,graphicx,url}

\usepackage{geometry}
\newcommand{\va}{\mathbf a}
\newcommand{\vb}{\mathbf b}

\newcommand{\vr}{\mathbf r}

\def\diag{\mathrm{diag}} %
\def\Span{\mathrm{span}} %
\newcommand{\divdif}{\mathord{\kern.43em{\vrule width.6pt height7pt depth-.28pt} \kern-.41em\Delta}}
\DeclareMathOperator*{\argmin}{arg\,min}

\newcommand*\vell{\ensuremath{\boldsymbol\ell}}

\usepackage{cellspace}
\title{Algorithms for the rational approximation of~matrix-valued~functions}
\author
{
	Ion  Victor Gosea\thanks{Max Planck Institute for Dynamics of Complex Technical Systems Magdeburg, 39106 Magdeburg, Germany,  \texttt{gosea@mpi-magdeburg.mpg.de}} 
	\and 
	Stefan G\"{u}ttel\thanks{The University of Manchester, Department of Mathematics, Manchester M13\,9PL, United Kingdon,  \texttt{stefan.guettel@manchester.ac.uk}. This  author has been supported by The Alan Turing Institute under the EPSRC grant EP/N510129/1.} 
}

\newcommand{\rev}[1]{#1}

\begin{document}

	\maketitle
	\begin{abstract}
		A selection of algorithms for the rational approximation of matrix-valued functions are discussed, including variants of the interpolatory AAA method, the RKFIT method based on approximate least squares fitting, vector fitting, and a method based on low-rank approximation of a block Loewner matrix. A new method, called the block-AAA algorithm, based on a generalized barycentric formula with matrix-valued weights is proposed. All algorithms are compared in terms of obtained approximation accuracy and runtime on a set of problems from model order reduction and nonlinear eigenvalue problems, including examples with noisy data. It is found that interpolation-based methods are typically cheaper to run, but they may suffer in the presence of noise for which approximation-based methods perform better.
	\end{abstract}
	

	\pagestyle{myheadings} \thispagestyle{plain} \markboth{ION VICTOR GOSEA AND  STEFAN G\"{U}TTEL}{MATRIX-VALUED RATIONAL APPROXIMATION}


	\section{Introduction} 
	Rational approximation is a powerful tool in applied science and engineering. To give just two examples, it is very commonly used for model order reduction~\cite{antoulas2005approximation,abg10,abg20} and the solution of nonlinear eigenvalue problems~\cite[Section~6]{GT17}. Recently, several new algorithms for the rational approximation and interpolation of scalar-valued functions have been proposed, including (in reverse chronological order) the AAA algorithm~\cite{nakatsukasa2018aaa}, the RKFIT algorithm~\cite{BeGu15g,BeGu15r}, vector fitting~\cite{GS99,BG06}, and methods based on the Loewner matrices~\cite{MA07}.  The aim of this paper is to explore extensions of these methods  for the purpose of approximating a matrix-valued (or block) function $F:\Lambda\to \mathbb{C}^{m\times n}$ on a discrete set $\Lambda$ in the complex plane. 
	
	The  paper contains two key contributions. Firstly, we  propose an extension of a barycentric formula to non-scalar weights and develop it into a new algorithm for the computation of rational approximants of matrix-valued functions. This algorithm, referred to as the block-AAA algorithm, generalizes the AAA algorithm in~\cite{nakatsukasa2018aaa}. Secondly, we perform  extensive numerical comparisons of several rational approximation algorithms  proposed by different communities (model reduction, numerical linear algebra, approximation theory). To the best of our knowledge, this is the first paper that offers a comprehensive study of these  algorithms. As part of our experiments, we identify potential problems of interpolation-based algorithms when the data samples are polluted by noise.
	
	The paper is structured as follows. Section~\ref{sec:scalar} contains a brief review of existing algorithms for scalar rational approximation. In Section~\ref{sec:represent}  we introduce new  representations of matrix-valued rational approximants and the new block-AAA algorithm. In Section~\ref{sec:alg} we briefly review other popular algorithms for rational approximation, including methods based on partial fractions, rational Krylov spaces, and state space representations. Section~\ref{sec:numer} contains a collection of numerical experiments and discussions, guided by six examples from different application areas. Section~\ref{sec:discussion} further discusses the performance of the different methods for each numerical example. We conclude in Section~\ref{sec:conclusion} with some ideas for future work.
	
	\section{Scalar rational approximation}\label{sec:scalar}
	In this section we summarize several algorithms for the rational approximation of a  scalar function $f:\Lambda\to\mathbb{C}$ sampled at a nonempty set of  points $\Lambda:=\{\lambda_1,\ldots,\lambda_\ell\}\subset \mathbb{C}$. This serves as an introduction of our notation, but also as a template for the block~variants in Sections~\ref{sec:represent} and~\ref{sec:alg}.
	
	\subsection{Adaptive Antoulas--Anderson (AAA) algorithm} \label{subsec:AAA_scalar}
	The AAA algorithm proposed in \cite{nakatsukasa2018aaa} is a practically robust and easy-to-use method for scalar rational interpolation. The degree~$d$ interpolant $r_d$ obtained after $d$ iterations of the AAA algorithm is of the form
	\begin{equation}\label{eq:bary}
	r_d(z) =  \frac{\displaystyle\sum_{k=0}^d \frac{w_k f_k}{z - z_k}}{\displaystyle\sum_{k=0}^d\frac{w_k}{z - z_k}},
	\end{equation}
	with nonzero barycentric weights $w_k\in\mathbb{C}$, pairwise distinct support points $z_k\in\mathbb{C}$, and function values $f_k = f(z_k)$. 
	A key ingredient of the AAA algorithm is its greedy choice of the support points, one per iteration $j=0,1,\ldots,d$, intertwined with the solution of a least squares problem to determine the barycentric weights $w_0^{(j)},\ldots,w_j^{(j)}$ at each iteration $j$. 
	The core AAA algorithm can be summarized as follows: 
	\medskip
	\begin{enumerate}
		\item Set $j=0$, $\Lambda^{(0)} := \Lambda$, and $r_{-1}:\equiv \ell^{-1}\sum_{i=1}^\ell f(\lambda_i)$.
		\item Select a  point $z_{j} \in \Lambda^{(j)}$ where $|f(z)-r_{j-1}(z)|$ is maximal, 
		with 
		\[
		r_{j-1}(z) :=  \frac{\displaystyle\sum_{k=0}^{j-1} \frac{w_k^{(j-1)} f_k}{z - z_k}}{\displaystyle\sum_{k=0}^{j-1}\frac{w_k^{(j-1)}}{z - z_k}} \quad \text{for} \  j\geq 1.
		\]
		\item If $|f(z_j) - r_{j-1}(z_j)|$ is small enough, return $r_{j-1}$ and stop.
		\item Set $f_j := f(z_j)$ and $\Lambda^{(j+1)} := \Lambda^{(j)}\setminus \{z_j\}$.
		\item Compute weights $w_0^{(j)},\ldots,w_j^{(j)}$ with $\sum_{k=0}^j |w_k^{(j)}|^2=1$ such that 
		\[
		\sum_{k=0}^j\frac{w_k^{(j)}}{z - z_k} f(z) \approx \sum_{k=0}^j \frac{w_k^{(j)} f_k}{z - z_k}
		\]
		is solved with least squares error over all $z\in\Lambda^{(j+1)}$.
		\item Set $j:=j+1$ and go to step 2.
	\end{enumerate}
	\medskip
	
	MATLAB implementations of the AAA algorithm can be found in \cite{nakatsukasa2018aaa} and the Chebfun package~\cite{chebfun} available at
	\smallskip
	\begin{center}
		\url{https://www.chebfun.org/}\\[1mm]
	\end{center}
	In \cite[Section~10]{nakatsukasa2018aaa} a non-interpolatory variant of the AAA algorithm is also mentioned, which is obtained by allowing the weights in the numerator and denominator of~\eqref{eq:bary} to be different.

	\subsection{Rational Krylov Fitting (RKFIT)} \label{subsec:RKFIT_scalar}
	
	The RKFIT method introduced in \cite{BeGu15g,BeGu15r} is based on iteratively relocating the poles of a rational function in the so-called RKFUN representation. An RKFUN is a triple $r_d\equiv (\underline{H_d},\underline{K_d}, \mathbf{c})$ with upper-Hessenberg matrices $\underline{H_d},\underline{K_d}\in\mathbb{C}^{(d+1)\times d}$ and a coefficient vector $\mathbf{c}\in\mathbb{C}^{d+1}$. For a given point $z\in\mathbb{C}$, the RKFUN $r_d(z)$ is evaluated as a dot-product $r_d(z):=\mathbf{n}(z)\cdot \mathbf{c}$, where $\mathbf{n}(z)$ is the unique left  nullvector of $z\underline{K_d} - \underline{H_d}$ normalized such that its first component is~$1$ (note that the unique existence of this nullvector is guaranteed by the fact that $z\underline{K_d} - \underline{H_d}$ is an unreduced upper-Hessenberg matrix). 
	
	The matrices $\underline{H_d}$ and $\underline{K_d}$ of the RKFUN representation satisfy  a rational Arnoldi decomposition $A V_{d+1} \underline{K_d} = V_{d+1} \underline{H_d}$ associated with a rational Krylov space  
	\[
	\mathcal{Q}_{d+1} (A,\mathbf{b}) := q(A)^{-1} \mathrm{span} \{ \mathbf{b}, A\mathbf{b}, \ldots, A^d \mathbf{b} \},
	\]
	with $A\in\mathbb{C}^{\ell\times \ell}$, $\mathbf{b}\in\mathbb{C}^\ell$, $V_{d+1}\in\mathbb{C}^{\ell\times (d+1)}$,  and a monic  polynomial $q\in\mathcal{P}_d$ such that $q(A)$ is invertible.  \rev{The columns of $V_{d+1}$ form an orthonormal basis of $\mathcal{Q}_{d+1} (A,\mathbf{b})$, generated by Ruhe's rational Krylov sequence (RKS) algorithm \cite{Ruhe84}.} It can be shown that the columns of $V_{d+1} = [\mathbf{v}_1,\ldots, \mathbf{v}_{d+1}]$ are all of the form $\mathbf{v}_{k+1} = p_k(A)q(A)^{-1}\mathbf{b}$ with polynomials $p_k\in\mathbb{P}_d$ for $k=0,1,\ldots,d$. In other words, the rational Arnoldi decomposition encodes a basis of scalar rational functions $r_k := p_k / q$ all sharing the same denominator $q$. Further, one can  show that the roots $\xi_j$ of $q(z) = \prod_{\xi_j \neq \infty} (z-\xi_j)$ correspond to the quotients of the subdiagonal elements of the upper-Hessenberg matrices $\underline{H_d}$ and $\underline{K_d}$, i.e., $\xi_j = h_{j+1,j}/k_{j+1,j}$ for $j=1,\ldots,d$. These quotients are referred to as the \emph{poles} of the rational Arnoldi decomposition \cite{BeGu15g}. 
	
	Given a matrix $F\in\mathbb{C}^{\ell\times \ell}$, the RKFIT algorithm attempts to identify an RKFUN $r_d\equiv (\underline{H_d},\underline{K_d}, \mathbf{c})$ such that $r_d(A)\mathbf{b} \approx F \mathbf{b}$ in the least squares sense.   Starting with an initial rational Arnoldi decomposition $A V_{d+1} \underline{K_d} = V_{d+1} \underline{H_d}$ having poles $\xi_1^{(0)},\ldots,\xi_d^{(0)}$, one RKFIT iteration consists of solving a least squares problem for finding a unit-norm vector $\widehat{\mathbf{v}}\in\Span (V_{d+1})$ such that
	\[
	\widehat{\mathbf{v}}	= \argmin_{\substack{\mathbf{v} = V_{d+1}\mathbf{c}\\ \|\mathbf{c}\|_2 = 1}} \| (I- V_{d+1}V_{d+1}^*) F\mathbf{v}\|_2,
	\]
	and to use the roots of $\widehat r = \widehat p/q$ associated with $\widehat{\mathbf{v}} = \widehat r(A)\mathbf{b}$ as the new poles $\xi_1^{(1)},\ldots, \hat \xi_d^{(1)}$ for the next iteration. This process is iterated until a convergence criterion is satisfied. 
	
	A MATLAB implementation of RKFIT is provided in the Rational Krylov Toolbox~\cite{BeGu14} available at
	\begin{center}
		\url{http://www.rktoolbox.org/}\\[2mm]
	\end{center}
	RKFIT can naturally deal with  scalar rational approximation  by choosing diagonal matrices $A = \diag(\lambda_1,\ldots, \lambda_\ell)$, $F=\diag(f_1,\ldots,f_\ell)$, and $\mathbf{b}=[1,1,\ldots,1]^T$. We  emphasize that RKFIT is a \emph{non-interpolatory} method.
	
	\subsection{Vector Fitting (VF)}\label{subsec:VF_scalar}
	
	The VF algorithm, originally proposed in \cite{GS99}, seeks to fit a rational function in partial fraction (pole--residue) form
	\begin{equation*}
	r_d(z) = \delta + \sum_{k=1}^{d} \frac{\gamma_k}{z-\xi_k}= \frac{p(z)}{q(z)}.
	\end{equation*}
	The iterative algorithm is initiated by choosing the degree $d$ of the rational approximant and an initial guess for the poles $\{\xi_1^{(0)}, \ldots, \xi_d^{(0)} \}$.
	At iteration $j=0,1,\ldots$ one determines parameters 
	$c_k^{(j)}$ and $d_k^{(j)}$ such that 
	\begin{equation*}
	\underbrace{\sum_{k=1}^{d} \frac{c_k^{(j)}}{z-\xi_k^{(j-1)}}+c_0^{(j)}}_{p^{(j)}(z)} \approx \Big{(} \underbrace{\sum_{k=1}^{d} \frac{d_k^{(j)}}{z-\xi_k^{(j-1)}}+1}_{q^{(j)}(z)} \Big{)} f(z)
	\end{equation*}
	is solved with least squares error over  all $z\in\Lambda$. Afterwards, the next set of poles $\{\xi_1^{(j)}, \ldots, \xi_d^{(j)} \}$ is computed as the roots of the polynomial $q^{(j)}(z)$ by solving a linear eigenvalue problem. The iteration continues until a convergence criterion is satisfied. Vector fitting is a \emph{non-interpolatory} method. 
	A MATLAB implementation of this method is available at\\[-2mm]
	\begin{center}
		\url{https://www.sintef.no/projectweb/vectorfitting/}
	\end{center}
	\smallskip
	\rev{It should be noted that an extension of VF to matrix-valued functions was also presented in \cite{DGB15}, together with a detailed analysis of numerical issues arising with particular test cases and connections to optimal $\mathcal{H}_2$ approximation.} 
	
	\subsection{Loewner framework (LF)}
	\label{subsec:Loewner_scalar}
	
	This  LF method uses the full data set in matrix format by forming (possibly very large) Loewner matrices. The original method, introduced in \cite{MA07}, is based on constructing a reduced-order rational function in  state-space representation. \rev{For additional details, we refer the reader to the tutorial papers \cite{ALI17} and \cite{morKarGA20a}, with the latter emphasizing applications on rational approximation, and also to the book \cite{abg20}}. In recent years, the Loewner framework has been extended to classes of mildly nonlinear systems, including bilinear systems \cite{AGI16}.
	
	
	Assuming that $\ell$, the number of sampling points, is an even positive integer, the first step in the Loewner framework is to partition  $\Lambda = \Lambda^L \cup \Lambda^R$ into two disjoint sets of the same cardinality. 
	The set $\Lambda^L$ contains the left points $\{x_1,\ldots,x_{\ell/2}\}$, while $\Lambda^R$ contains the right points $\{y_1,\ldots,y_{\ell/2}\}$. Similarly, the set of scalar samples (evaluations of the function $f$ on $\Lambda$) is partitioned into two sets.
	
	One then defines matrices $\mathbb{L} \in \mathbb{C}^{(\ell/2) \times (\ell/2)}$ (the Loewner matrix) and also $\mathbb{L}^s \in \mathbb{C}^{(\ell/2) \times (\ell/2)}$ (the shifted Loewner matrix) with entries
	\[
	\mathbb{L}_{i,j} = \frac{ f(x_i) - f(y_j)}{x_i-y_j} \quad \text{and} \quad \mathbb{L}^s_{i,j} = \frac{x_i f(x_i) - y_j  f(y_j)}{x_i-y_j},
	\]
	respectively.  Additionally, one defines vectors $V \in \mathbb{C}^{ (\ell/2) \times 1}, \ W \in \mathbb{C}^{1 \times (\ell/2)}$  as
	\begin{align}\label{VW_mat_sca}
	V &= \left[  f(x_1),  f(x_2), \ldots, f(x_{\ell/2})  \right]^T, \ \ \ 
	W = \left[ f(y_1) , f(y_2) , \ldots , f(y_{\ell/2}) \right].
	\end{align}
	The next step is to compute a rank-$d$ truncated  singular value decomposition of the Loewner matrix $\mathbb{L}\approx X S Z^*$ with $X, Z \in \mathbb{C}^{(\ell/2) \times d}$, and  $S \in \mathbb{C}^{d \times d}$. 
	Finally, by means of projecting with the matrices $X$ and $Z$, the fitted rational function is  
	\begin{equation}\label{rd_Loewner}
	r_d(z) = W Z \left( X^* (\mathbb{L}^s - z \mathbb{L} )  Z \right)^{-1} X^* V.
	\end{equation}
	This is a subdiagonal rational function of type $(d-1,d)$ whose $d$ poles are given by the generalized eigenvalues of the matrix pair $(X^* \mathbb{L}^s Z, X^* \mathbb{L} Z)$. 
	
	\rev{It is to be noted that we will be using a slightly different implementation than that presented in \cite{MA07,ALI17}. More precisely, matrices $X$ and $Z$ in (\ref{rd_Loewner}) are computed solely from the SVD of $\mathbb{L}$, and not from the SVDs of the two "augmented Loewner matrices" $\left[\begin{matrix}
		\mathbb{L} \\ \mathbb{L}_s
		\end{matrix} \right] \in \mathbb{C}^{\ell \times (\ell/2)}$ and $\left[\begin{matrix}
		\mathbb{L} & \mathbb{L}_s
		\end{matrix} \right] \in \mathbb{C}^{(\ell/2) \times \ell}$, as originally proposed in \cite{MA07}.}
	
	%
	%

	\section{Matrix-valued barycentric forms and block-AAA}\label{sec:represent}
	
	The simplest matrix-valued barycentric form is obtained from \eqref{eq:bary} by replacing the function values $f_k$ with matrices $F_k := F(z_k)$:
	\begin{equation}
	\tag{bary-A}
	R_d(z) = \frac{\displaystyle\sum_{k=0}^d \frac{w_k F_k}{z - z_k}}{\displaystyle\sum_{k=0}^d\frac{w_k}{z - z_k}}.\label{eq:bb1}
	\end{equation}
	Provided that all weights $w_k$ are nonzero, the function $R_d$ interpolates the function~$F$ at all the support points~$z_k$. Each $(i,j)$ entry of $R_d$ is a rational function of the form $p_{ij}(z)/q(z)$, where $p_{ij},q\in\mathcal{P}_d$ are polynomials of degree $d$. Note that all these entries share the same \emph{scalar} denominator~$q$.

	A slight modification of \eqref{eq:bb1} yields a new matrix-valued barycentric formula  
	\begin{equation}
	\tag{bary-B}
	R_d(z) = \left( \sum_{k=0}^d \frac{W_k}{z-z_k}\right)^{-1} \left( \sum_{k=0}^d \frac{W_k F(z_k)}{z-z_k}\right) 
	\label{eq:bb2}
	\end{equation}
	with weight matrices $W_k\in\mathbb{C}^{m\times m}$. If all these $W_k$ are nonsingular, \rev{which needs to be assumed}, $R_d$ interpolates $F$ at all the support points $z_k$. 
	Given a set of support points $z_0,\ldots,z_k$ different from any point in $\Lambda$, a linearized version of the approximation problem $R_d(\lambda_i)\approx F(\lambda_i)$ is 
	\[
	\sum_{i=1}^\ell \left\|  \left( \sum_{k=0}^d \frac{W_k}{\lambda_i -z_k}\right) F(\lambda_i) - \sum_{k=0}^d \frac{W_k F(z_k)}{\lambda_i-z_k}   \right\|_F^2 \to \min_{W_k}. 
	\]
	The weight matrices can now be obtained from a trailing left-singular block vector $[ W_0,\ldots,W_d ]$ of unit  norm such that 
	\[
	[ W_0,\ldots,W_d ] 
	\begin{bmatrix}
	\frac{F(\lambda_1) - F(z_0)}{\lambda_1-z_0} & \cdots & \frac{F(\lambda_\ell) - F(z_0)}{\lambda_\ell-z_0} \\
	\vdots & & \vdots \\
	\frac{F(\lambda_1) - F(z_d)}{\lambda_1-z_d} & \cdots & \frac{F(\lambda_\ell) - F(z_d)}{\lambda_\ell-z_d}
	\end{bmatrix} = :\mathbb{W} \mathbb{L}
	\]
	has smallest possible \rev{Frobenius} norm. Note that the matrix  $\mathbb{L}$ is a block Loewner matrix.\footnote{In the Loewner framework discussed in Section~\ref{subsec:Loewner_scalar}, the Loewner matrix entries 
		\[
		\mathbb{L}_{ij} = \frac{\mathbf{v}_i^* \mathbf{r}_j - \vell_i^* \mathbf{w}_j}{\mu_i-\lambda_j}
		\]
		are defined in terms of left and right directions $(\mu_i,\vell_i,\mathbf{v}_i)$ and  $(\lambda_j,\mathbf{r}_j,\mathbf{w}_j)$, but here we have a special case where the left and right directions $\vell_i$ and $\mathbf{r}_j$ are chosen as the unit vectors.} 
	It is now easy to derive an AAA-like algorithm based on \eqref{eq:bb2}. It is shown in \rev{Figure~\ref{alg:baaa}} and referred to as the \emph{block-AAA algorithm}. 
	
	\begin{figure}
		\noindent\fbox{%
			\parbox{.97\textwidth}{%
				\textbf{Block-AAA algorithm}\\[1mm]
				\textbf{Inputs:} \,\,Discrete set $\Lambda\subset \mathbb{C}$ with $\ell$ points, function $F$, error tolerance $\varepsilon>0$\\[-0mm]
				\textbf{Output:} Rational approximant $R_{j-1}$ in the form \eqref{eq:bb2}\\[-3mm]
				\begin{enumerate}
					\item Set $j=0$, $\Lambda^{(0)} := \Lambda$, and $R_{-1}:\equiv \ell^{-1}\sum_{i=1}^\ell F(\lambda_i)$.
					\item Select a  point $z_{j} \in \Lambda^{(j)}$ where $\|F(z)-R_{j-1}(z)\|_F$ is maximal, 
					with 
					\[
					R_{j-1}(z) :=  \left(\displaystyle\sum_{k=0}^{j-1}\frac{W_k^{(j-1)}}{z - z_k}\right)^{-1} \left(\displaystyle\sum_{k=0}^{j-1} \frac{W_k^{(j-1)} F_k}{z - z_k}\right) \quad \text{for} \  j\geq 1.
					\]
					\item If $\|F(z_j) - R_{j-1}(z_j)\|_F\leq \varepsilon$, return $R_{j-1}$ and stop.
					\item Set $F_j := F(z_j)$ and $\Lambda^{(j+1)} := \Lambda^{(j)}\setminus \{z_j\}$.
					\item Find matrix $[W_0^{(j)},\ldots,W_j^{(j)}]$ of unit Frobenius norm such that 
					\[
					\sum_{k=0}^j\frac{W_k^{(j)}}{z - z_k} F(z) \approx \sum_{k=0}^j \frac{W_k^{(j)} F_k}{z - z_k}
					\]
					is solved with least squares error over all $z\in\Lambda^{(j+1)}$.
					\item Set $j:=j+1$ and go to step 2.
				\end{enumerate}
			}
		}
		\caption{Pseudocode of the block-AAA algorithm\label{alg:baaa}}
	\end{figure}

	A MATLAB implementation of the block-AAA algorithm is available at
	\smallskip
	\begin{center}
		\url{https://github.com/nla-group/block_aaa} \\[1mm]
	\end{center}
	As it will become clearer from the discussions in Section~\ref{sec:alg}, the block-AAA algorithm is different from the set-valued AAA \cite{LPVM18} and the fast-AAA \cite{Hoch17} algorithms: the entries of a block-AAA approximant $R_d$ \emph{do not} share a common scalar denominator of degree~$d$. Indeed, a block-AAA approximant $R_d$ of order $d$ can have a larger number of up to $dm$ poles, i.e., its McMillan degree can be as high as $dm$. For this reason we refer to the \emph{order} of a block-AAA approximant instead of its degree.

	Finally, we mention for completeness two other barycentric formulas which might be of interest elsewhere. 
	One is given~as 
	\begin{equation}
	\tag{bary-C}
	R_d(z) = \left( \sum_{k=0}^d \frac{D_k}{z-z_k}\right)^{-1} \left( \sum_{k=0}^d \frac{C_k}{z-z_k}\right) 
	\label{eq:bb3}
	\end{equation}
	with the matrices $C_k\in\mathbb{C}^{m\times n}$ and $D_k\in\mathbb{C}^{m\times m}$ chosen such that $R_\ell(\lambda_i)\approx F(\lambda_i)$. We  assume that the $z_k$ and $\lambda_i$ are pairwise distinct and consider the linearized problem
	\[
	\sum_{i=1}^\ell \left\|  \left( \sum_{k=0}^d \frac{D_k}{\lambda_i -z_k}\right) F(\lambda_i) - \sum_{k=0}^d \frac{C_k}{\lambda_i-z_k}   \right\|_F^2 \to \min_{C_k, D_k}. 
	\]
	This least squares problem is solved by a left-singular block vector given by $[ C_0,\ldots,C_d,D_0,\ldots,D_d ]$ of unit \rev{Frobenius} norm such that 
	\[
	[ C_0,\ldots,C_d,D_0,\ldots,D_d ] 
	\begin{bmatrix}
	\frac{-I}{\lambda_1-z_0} & \cdots & \frac{-I}{\lambda_\ell-z_0} \\
	\vdots & & \vdots \\
	\frac{-I}{\lambda_1-z_d} & \cdots & \frac{-I}{\lambda_\ell-z_d} \\[2mm]
	\frac{F(\lambda_1)}{\lambda_1-z_0} & \cdots & \frac{F(\lambda_\ell)}{\lambda_\ell-z_0} \\
	\vdots & & \vdots \\
	\frac{F(\lambda_1)}{\lambda_1-z_d} & \cdots & \frac{F(\lambda_\ell)}{\lambda_\ell-z_d}
	\end{bmatrix}
	\]
	has smallest possible \rev{Frobenius} norm.  Provided that the matrices $D_k$ are nonsingular, \rev{which needs to be assumed}, we have $R_d(z_k) = D_k^{-1} C_k$ and in general the approximation is \emph{non-interpolatory}. Such approximants might be useful for the solution of nonlinear eigenvalue problems; see also the appendix.
	
	It is also possible to have matrix-valued "support points" $Z_k\in\mathbb{C}^{m\times m}$, resulting~in 
	\begin{equation}
	\tag{bary-D}
	R_d(z) = \left( \sum_{k=0}^d (z I - Z_k)^{-1} D_k\right)^{-1} 
	\left( \sum_{k=0}^d (z I - Z_k)^{-1} C_k\right).
	\label{eq:bb4}
	\end{equation}
	This leads to tangential interpolation determined by the eigenvalues and eigenvectors of the matrices $Z_k$. We do not further explore this form here. \rev{Note that tangential interpolation was previously used for model reduction of MIMO systems in \cite{GVD04}.}

	\section{Other block methods}\label{sec:alg}
	
	In this section we list methods for the rational approximation of matrix-valued functions using the barycentric formula \eqref{eq:bb1}, as well as methods based on other  representations of their approximants. These methods  will be compared numerically in  Sections~\ref{sec:numer} and~\ref{sec:discussion}. 
	
	\smallskip
	
	\noindent\textbf{Set-valued AAA (fast-AAA):} The set-valued AAA algorithm presented in~\cite{LPVM18}, and similarly the fast-AAA algorithm in \cite{Hoch17}, applies the standard AAA algorithm to each component of $F$ using common weights and support points for all of them, thereby effectively producing a barycentric interpolant in the form~\eqref{eq:bb1}. 
	
	\smallskip
	
	\noindent\textbf{Surrogate AAA:} This method presented in~\cite{EG17} applies the standard AAA algorithm to a scalar surrogate function $f(z) := \va^T F(z) \vb$, with vectors $\va,\vb$ chosen at random. The resulting block interpolant is of the form~\eqref{eq:bb1}, obtained by replacing the scalar values $f_k = f(z_k)$ by the matrices $F_k = F(z_k)$. 
	
	\smallskip

	\noindent\textbf{RKFIT (block):} The RKFIT algorithm has been generalized in \cite{BeGu15r} to the problem of finding a family $r_d^{[k]}(z) = p^{[k]}(z)/q(z)$ of scalar rational functions such that $r_d^{[k]}(z) \approx f^{[k]}(z)$ in the least squares sense over all $z\in\Lambda$. All computed functions $r_d^{[k]}$ are represented in the RKFUN format and share the common scalar denominator $q$ of degree $d$. In our setting, we simply associate with each $(i,j)$ entry of $F$ a corresponding function  $f^{[k]}(z)$ and run RKFIT on that family. 
	
	\smallskip

	\noindent\textbf{Matrix Fitting:} 
	The extension of vector fitting to matrix-valued data has been implemented in the MATLAB Matrix Fitting Toolbox   available at
	\smallskip
	\begin{center}
		\url{https://www.sintef.no/projectweb/vectorfitting/downloads/matrix-fitting-toolbox/} \\[1mm]
	\end{center} 
	As described in the manual \cite{MatVF} associated with the toolbox, the primary intention of this  software is the rational multi-port modeling of data in the frequency domain. In particular, the toolbox deals with fitting so-called admittance values (known also as $Y$ parameters) and scattering values (known also as $S$ parameters). 
	
	The fitted rational function is provided in the pole--residue form
	\begin{equation*}
	R_d(z) = D + \sum_{k=1}^{d} \frac{C_k}{z-\xi_k},
	\end{equation*}
	or as an equivalent state space model
	\begin{equation*}
	R_d(z) = D + C (z I- A)^{-1} B,
	\end{equation*}
	where $D,C_k \in \mathbb{C}^{n \times n}$ for $k = 1,\ldots,d$, $C \in \mathbb{C}^{n \times n d}$, $A \in \mathbb{C}^{n d \times n d}$, and $B \in \mathbb{C}^{n d \times n}$.
	For both representations, stable poles can be enforced. Note that the matrix elements of $R_d$ share a common set of $d$ poles. 
	
	The above mentioned matrix fitting implementation  \cite{MatVF} works with the elements of the upper triangular part of the matrix samples $F_i := F(\lambda_i)$ ($i=1,\ldots,\ell$), i.e., the entries $[F_i]_{a,b}$ for $a \leq b$, which are stacked into a single vector. Hence, this implementation can be used only for symmetric matrix-valued functions. 
	The elements of the vector are fitted using an implementation of the ``relaxed version'' of vector fitting~\cite{BG06}, provided in the function \texttt{vectfit3.m} of the toolbox. The main driver function \texttt{VFdriver.m} accepts as input arguments sample points $\Lambda = \overline{\Lambda} = \{\lambda_1,\ldots,\lambda_\ell\}$ on the imaginary axis in complex conjugate pairs (see Section~3.1, page~7 of the toolbox manual). In our context of fitting on arbitrary discrete sets $\Lambda$ in the complex plane, this requirement might represent an inconvenience. Also, the available driver seems to be  intended to be run for even orders~$d$ only, although this is not an inherent restriction of vector fitting. 
	
	\smallskip
	
	\noindent\textbf{Block Loewner framework:} The extension of the Loewner framework to matrix-valued functions is covered in \cite{MA07}. Therein, the case of directionally-sampled matrix data is treated in detail. Introducing left directional vectors $\vell_i \in \mathbb{C}^m$ and right directional vectors $\vr_j \in \mathbb{C}^n$, one defines a Loewner matrix  $\mathbb{L} \in \mathbb{C}^{(\ell/2) \times (\ell/2)}$ and a shifted Loewner matrix as follows:
	\begin{align*}
	\mathbb{L}_{i,j} &= \frac{\vell_i^* F(x_i) \vr_j-\vell_i^* F(y_j) \vr_j}{x_i-y_j}, \\ \mathbb{L}^s_{i,j} &= \frac{x_i \vell_i^* F(x_i) \vr_j- y_j \vell_i^* F(y_j) \vr_j}{x_i-y_j}.
	\end{align*}
	The definition of the matrices $V \in \mathbb{C}^{ (\ell/2) \times n}, \ W \in \mathbb{C}^{m \times (\ell/2)}$ is also different from those in the scalar case in (\ref{VW_mat_sca}), namely the directional vectors appear in the form  
	\begin{align*}
	V^* &= \left[ \begin{matrix} F(x_1)^* \vell_1  , F(x_2)^* \vell_2 , \ldots , F(x_{\ell/2})^* \vell_{\ell/2} \end{matrix} \right], \\
	W &= \left[ \begin{matrix} F(y_1) \vr_1 , F(y_2) \vr_2 , \ldots , F(y_{\ell/2}) \vr_{\ell/2} \end{matrix} \right].
	\end{align*}
	The rational matrix-valued function $R_d$ is then computed exactly as in the scalar case  (\ref{rd_Loewner}) in the form 
	\begin{equation*}
	R_d(z) = W Z \left( X^* (\mathbb{L}^s - z \mathbb{L} )  Z \right)^{-1} X^* V,
	\end{equation*}
	where the truncated singular vector matrices $X,Z$ are the same as in Section \ref{subsec:Loewner_scalar}.

	\section{Numerical comparisons}\label{sec:numer}
	
	In all our tests we assume that $\Lambda = \{ \lambda_1,\ldots,\lambda_\ell\}$ and  all the samples $F(\lambda_i)$ are given. 
	We will compare the six algorithms introduced in the previous sections, using MATLAB implementations provided by their authors whenever available.  
	In order to compare the algorithms with respect to their approximation performance, we use the root mean squared error
	\begin{equation*}
	\text{\texttt{RMSE}} = \left( \frac{1}{\ell} \sum_{i=1}^\ell \| F(\lambda_i) - R_d(\lambda_i) \|_F^2\right)^{1/2}. 
	\end{equation*}
	As before, $d$ generally refers to the order of the rational approximant $R_d$, but it is the same as the degree for all algorithms except block-AAA. The RMSE values are reported for various examples and different different orders~$d$ in Tables~\ref{table_example2}, \ref{table_example3}, \ref{table_example4}, \ref{table_example5}, and~\ref{table_example6}. 
	
	We also report timings for all algorithms under consideration. Although MATLAB timings might not always be  reliable, the differences in the runtimes are usually large enough to give some indication of algorithmic complexity. In order to reduce random fluctuations in the timings, we run each algorithm for 20 times and compute the average execution time for a single run. 
	All experiments were performed on a desktop computer with \rev{16~GB RAM and an Intel(R) Core(TM) i7-10510U CPU running at 1.80~GHz 2.30~GHz.}

	\subsection{Two toy examples}\label{subsec:acad}
	We first consider a rational function given by 
	$$
	F(z) = \left[ \begin{matrix} \frac{2}{z+1} & \frac{3-z}{z^2+z-5} \\[1mm]
	\frac{3-z}{z^2+z-5} & \frac{2+z^2}{z^3 + 3z^2-1} \end{matrix} \right].
	$$
	The entries of $F$ can be written with a  common scalar denominator of degree $d=6$, e.g., $q(z) = (z+1)(z^2+z-5)(z^3 + 3z^2-1)$.  We use $\ell = 100$ logarithmically-spaced samples of $F$ on the imaginary axis in the interval $[1,100]\mathrm{i}$. 
	As in addition the samples of $F$ are symmetric matrices, this example is suitable for the  Matrix Fitting Toolbox~\cite{MatVF}. 
	
\begin{figure}[h]
		\hspace*{-0mm}\includegraphics[scale=.33]{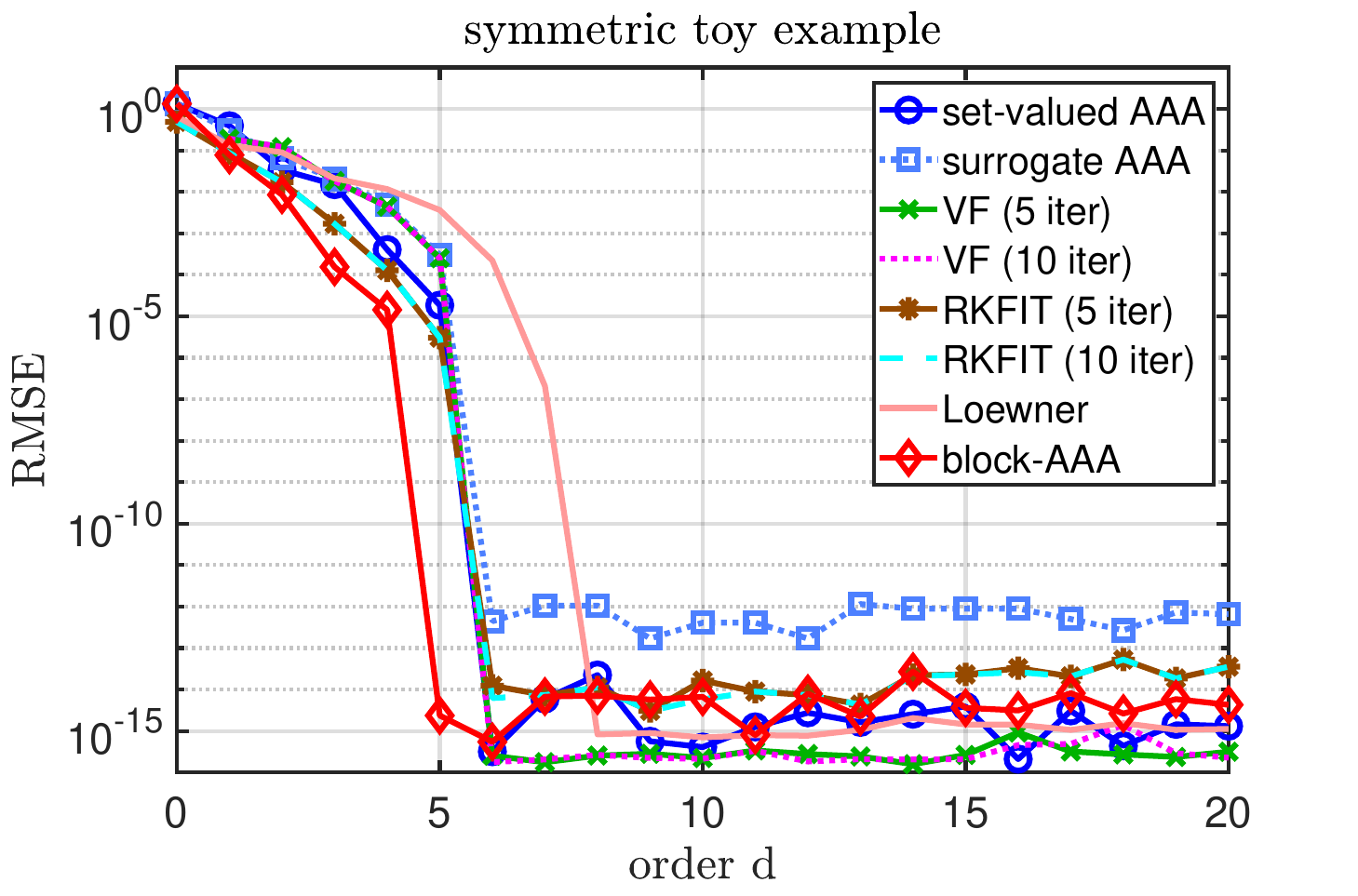}\hspace{-3mm}\includegraphics[scale=.33]{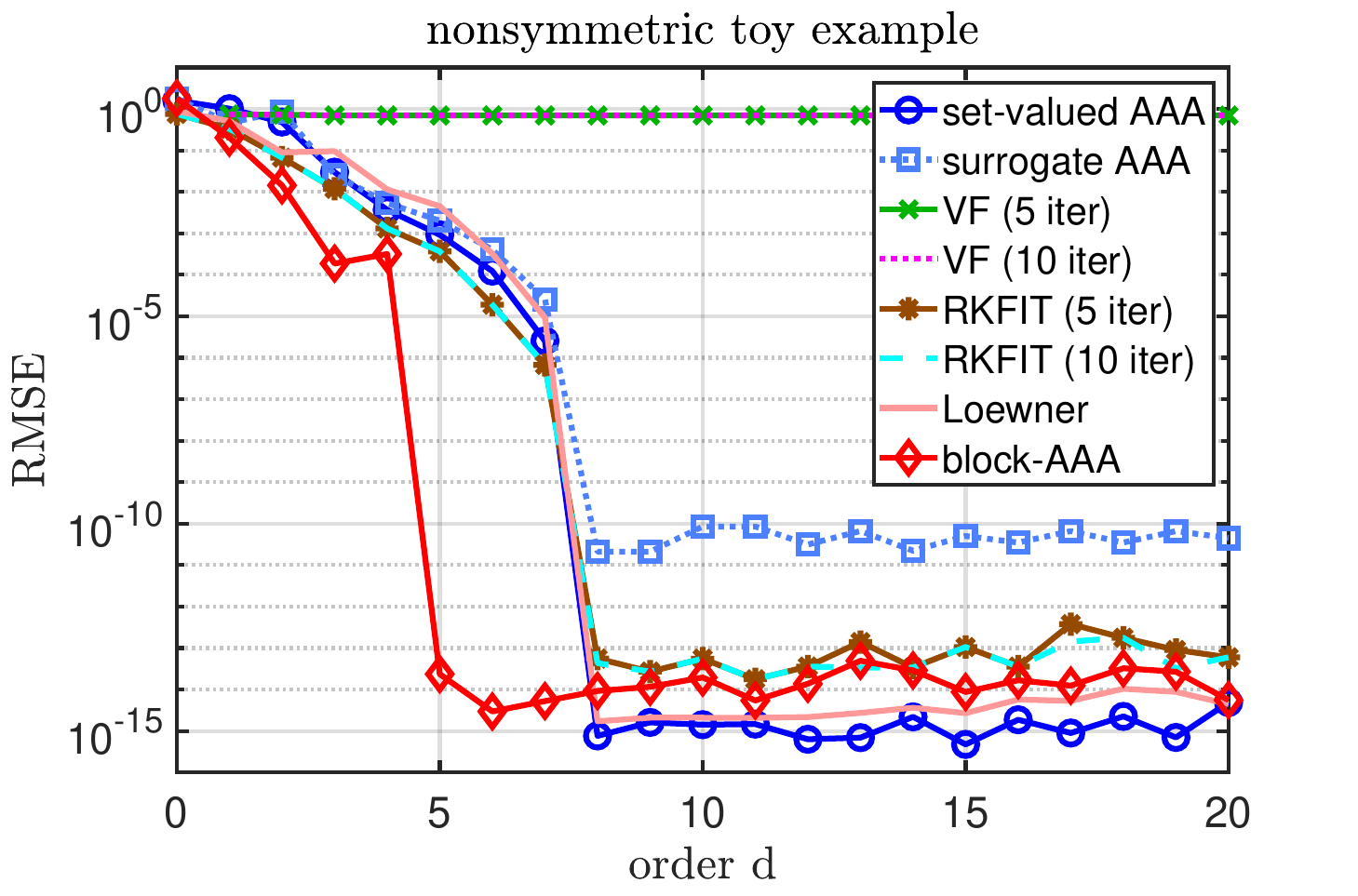}
		\vspace{-6mm}
		\caption{Accuracy comparison for the two toy examples in Section~\ref{subsec:acad}. Left: approximating a symmetric $2\times 2$ matrix-valued function. Right: modified non-symmetric case.}
		\label{fig:example_mftex1}
		\vspace{-0mm}
	\end{figure}

	The RMSE values obtained with all discussed algorithms for varying orders $d=0,1,\ldots, 20$ are shown in Figure~\ref{fig:example_mftex1} (left). 
	Note that all methods eventually recover the function accurately. The Loewner approach  requires a slightly higher order than the other methods, i.e., $d_1 = 8$. The block-AAA algorithm correctly identifies $F$ when the order is $d_2 = 5$.

	Next, we modify the (1,2) entry of the function $F$ by replacing its denominator by  $z^2+z+5$. The entries of $F$ can then be written with a common denominator of degree  $d=8$, e.g., $q(z) = (z+1)(z^2+z-5)(z^2+z+5)(z^3 + 3z^2-1)$.
	The RMSE values are shown in Figure~\ref{fig:example_mftex1} (right). Now the VF approach fails as the modified function~$F$ is no longer symmetric. The block-AAA algorithm again identifies $F$ correctly with an order of  $d_1 = 5$, while all other methods require a degree of $d_2 = 8$ as expected.
	
	\subsection{An example from the MF Toolbox}\label{subsec:matfit} 
	
\begin{figure}[h]
		\includegraphics[scale=.33]{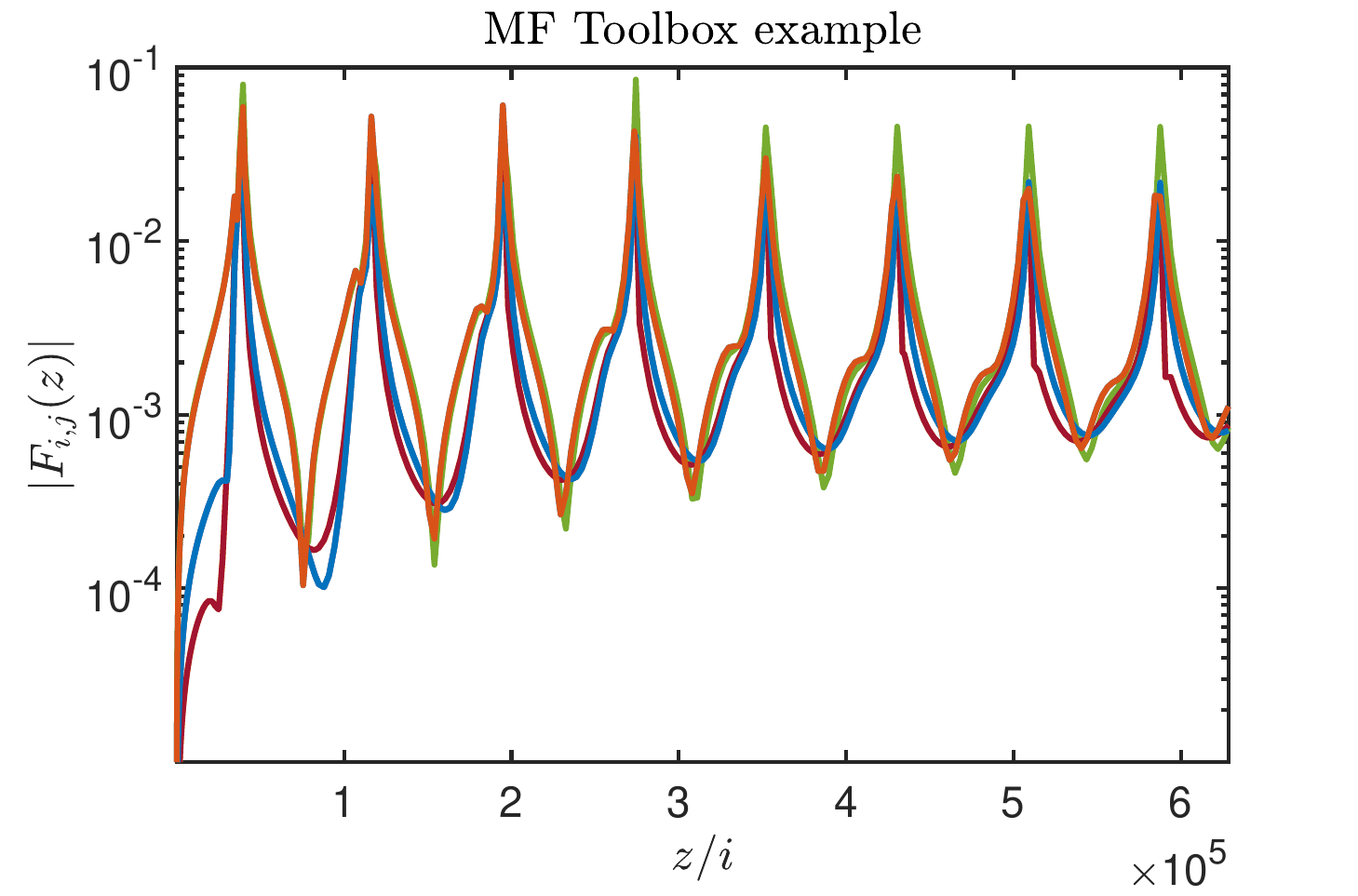}\hspace*{-3mm}\includegraphics[scale=.33]{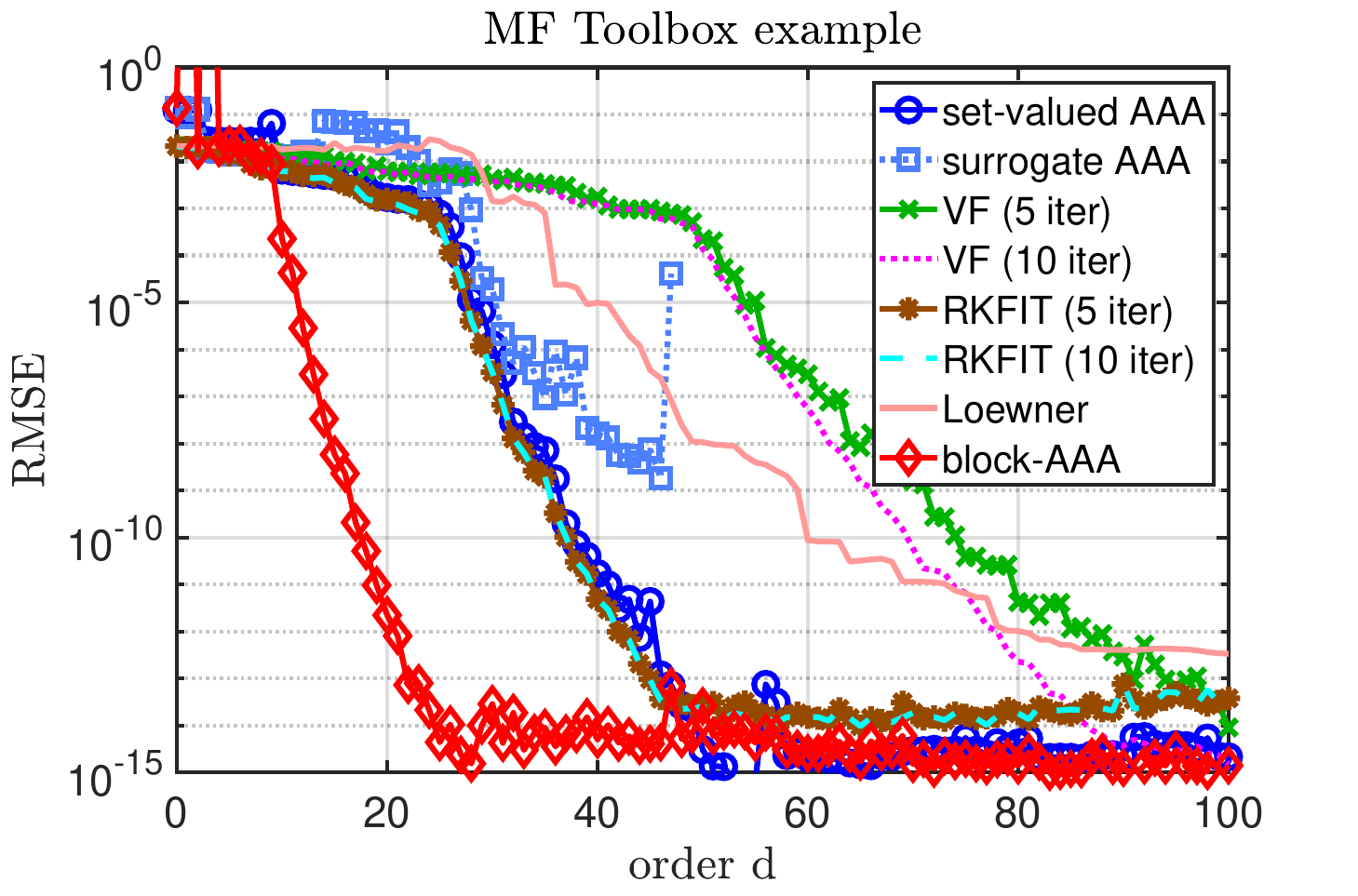}
		\vspace{-6mm}
		\caption{Example with admittance matrices from the MF Toolbox, described in Section~\ref{subsec:matfit}. Left: entries of the $3\times 3$  matrix $F(z)$ over the frequency range. Right: accuracy performance.}
		\label{fig:example_mftex22}
		\vspace{-0mm}
	\end{figure}

	We now choose an example from the Matrix Fitting Toolbox~\cite{MatVF}. We use the file \texttt{ex2\_Y.m} provided therein, containing 300 samples of $3 \times 3$ admittance matrices for a three-conductor overhead line.
	The sampling points are all on the imaginary axis within the interval $62.8319 \cdot [1,10^4]$i. Figure~\ref{fig:example_mftex22} (left) depicts the magnitude of the nine matrix entries over the sampled frequency range.

	Figure~\ref{fig:example_mftex22} (right) shows the RMSEs achieved by each of the tested algorithms for varying orders $d=0,1,\ldots,50$. For two particular orders, namely $d=10$ and $d=20$, we report the numerical RMSE values and the corresponding timings ("runtime") in Table~\ref{table_example2}. \rev{Additionally, in the same table we display the evaluation times ("evaltime") for each method; this represents the time needed for evaluating a rational function at 1000 values, logarithmically-spaced in the original sampling interval.}
	
	\rev{For another test, reported in Tables~\ref{table_example22} and \ref{table_example222}, we choose four target RMSE values $\epsilon_1 = 10^{-3}, \ \epsilon_2 = 10^{-6}, \ \epsilon_3 = 10^{-9}$ and $\epsilon_4 = 10^{-12}$. We then we report the orders and runtimes needed by each of the eight methods to achieve these RMSEs. In two cases, the surrogate AAA method fails to achieve the target RMSE and this is denoted by an ``--''.}

	\begin{table}[h] 
		\rev{\small \caption{Selected RMSE values and timings for all tested algorithms --- MF Toolbox  example}
			\label{table_example2}
			\begin{center}
				\begin{tabular}{|l|c|c|c|c|c|l|}\hline
					& \multicolumn{2}{c|}{RMSE} & \multicolumn{2}{c|}{Runtime (ms)} & \multicolumn{2}{c|}{Evaltime (ms)}\\
					& $d=10$ & $d=20$ & $d=10$ & $d=20$ & $d=10$ & $d=20$  \\\hline
					Set-valued AAA  & $6.522 \cdot 10^{-3}$   & $1.707 \cdot 10^{-3}$ &  $11.5$   & $24.2$ &  $18.5$ & $27.6$   \\
					Surrogate AAA  & $1.278 \cdot 10^{-2}$    & $2.501 \cdot 10^{-2}$ &  $9.7$   & $15.9$  & $15.6$ & $23.6$  \\
					VF (5 iter) &  $1.374 \cdot 10^{-2}$ &  $6.175 \cdot 10^{-3}$ &   $27.7$   & $46.6$  & $21.8$ & $39.9$  \\
					VF (10 iter) &  $1.139 \cdot 10^{-2}$   & $6.146 \cdot 10^{-3}$ &  $41.9$   & $70.6$   & $21.9$ & $40.8$   \\
					RKFIT (5 iter) &  $6.619 \cdot 10^{-3}$   & $1.568 \cdot 10^{-3}$ &  $68.1$   & $134.2$  & $51.3$ & $122.9$   \\
					RKFIT (10 iter) &  $6.246 \cdot 10^{-3}$   & $1.233 \cdot 10^{-3}$ &  $110.9$   & $225.2$  & $54.7$ & $119.8$   \\
					Loewner &  $1.795 \cdot 10^{-2}$   & $1.675 \cdot 10^{-2}$ &  $18.6$   & $19.5$  & $15.8$ & $30.5$   \\
					block-AAA &  $2.255 \cdot 10^{-4}$   & $2.215 \cdot 10^{-12}$  &  $98.1$   & $207.0$  & $60.4$  &  $79.5$ \\
					\hline
				\end{tabular}
		\end{center}}
	\end{table}
	
	\normalsize
	
	\begin{table}[h] 
		\rev{	\caption{The orders for 4 different RMSE values $\epsilon_k, \ k = 1,4$ for all algorithms --- MF Toolbox  example}
			\label{table_example22}
			\begin{center}
				\begin{tabular}{|l|c|c|c|c|c|l|}\hline
					& \multicolumn{4}{c|}{Orders} \\
					& $\epsilon_1 = 10^{-3}$ & $\epsilon_2 = 10^{-6}$ & $\epsilon_3 = 10^{-9}$ & $\epsilon_4 = 10^{-12}$   \\\hline
					Set-valued AAA & $25$ & $31$ & $37$ & $45$  \\
					Surrogate AAA  & $28$ & $34$ & $-$ & $-$   \\
					VF (5 iter) &   $41$ & $56$ &  $71$ &  $86$   \\
					VF (10 iter) & $41$ & $55$ & $66$ & $77$     \\
					RKFIT (5 iter) & $23$ & $30$ & $36$ & $42$   \\
					RKFIT (10 iter) &  $22$ & $30$ & $36$ & $42$   \\
					Loewner &  $35$ & $45$ & $60$ & $81$   \\
					block-AAA &  $10$ & $13$ & $17$ & $21$  \\
					\hline
				\end{tabular}
		\end{center}}
	\end{table}
	
	\begin{table}[h] 
		\rev{	\caption{The timings for 4 different RMSE values $\epsilon_k, \ k = 1,4$ for all algorithms --- MF Toolbox  example}
			\label{table_example222}
			\begin{center}
				\begin{tabular}{|l|c|c|c|c|c|l|}\hline
					& \multicolumn{4}{c|}{Runtime (ms)} \\
					& $\epsilon_1 = 10^{-3}$ & $\epsilon_2 = 10^{-6}$ & $\epsilon_3 = 10^{-9}$ & $\epsilon_4 = 10^{-12}$
					\\\hline
					Set-valued AAA & $23.2$ & $29.2$ & $37.7$ & $49.1$  \\
					Surrogate AAA  & $19.5$ & $25.3$ & $-$ & $-$   \\
					VF (5 iter) &  $74.5$ & $110.2$ & $127.2$ & $178.4$   \\
					VF (10 iter) & $135.4$ & $169.9$ & $209.9$ & $247.1$     \\
					RKFIT (5 iter) & $154.9$ & $221.5$ & $287.7$ & $366.6$   \\
					RKFIT (10 iter) &  $229.0$ & $356.5$ & $465.1$ & $559.6$   \\
					Loewner &  $25.2$ & $42.2$ & $56.4$ & $120.2$   \\
					block-AAA &  $78.3$ & $106.5$ & $142.5$ & $197.2$  \\
					\hline
				\end{tabular}
		\end{center}}
	\end{table}


	\subsection{A first model order reduction example} \label{subsec:CD}
	
	We now consider the \texttt{CD player} example from the SLICOT benchmark collection  \cite{bench02}. The mechanism to be modeled
	is a swing arm on which a lens is mounted by means of two horizontal leaf springs. The challenge is to find a low-cost controller that can make the servo-system faster and less sensitive to external shocks. The LTI system that models the dynamics has 60 vibration modes, hence the dimension $n=120$. Additionally, there are $m=2$ inputs and $p=2$ outputs. See \cite{CVD05} for more details.
	
\begin{figure}[h]
		\includegraphics[scale=.33]{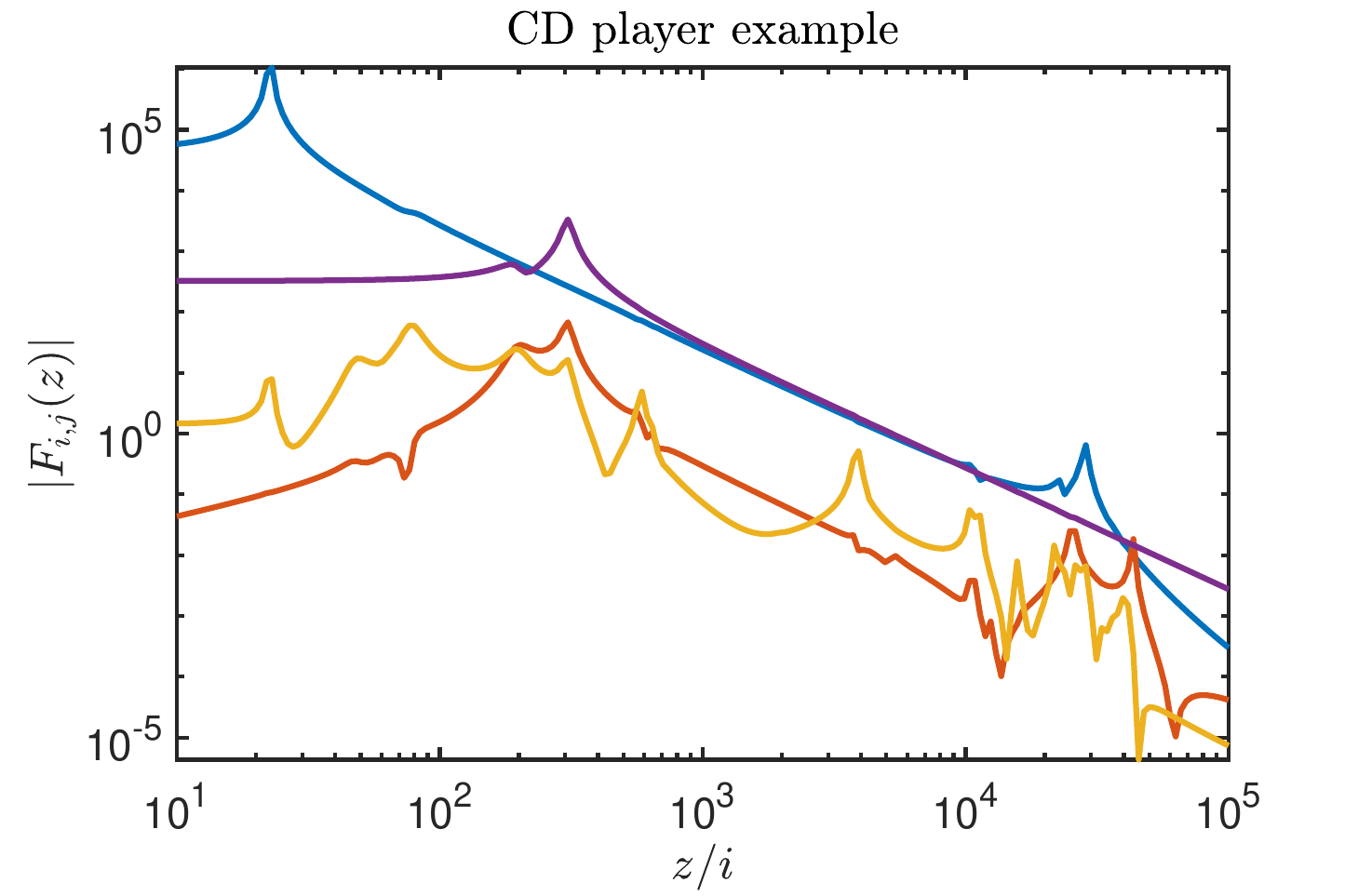}\hspace*{-3mm}\includegraphics[scale=.33]{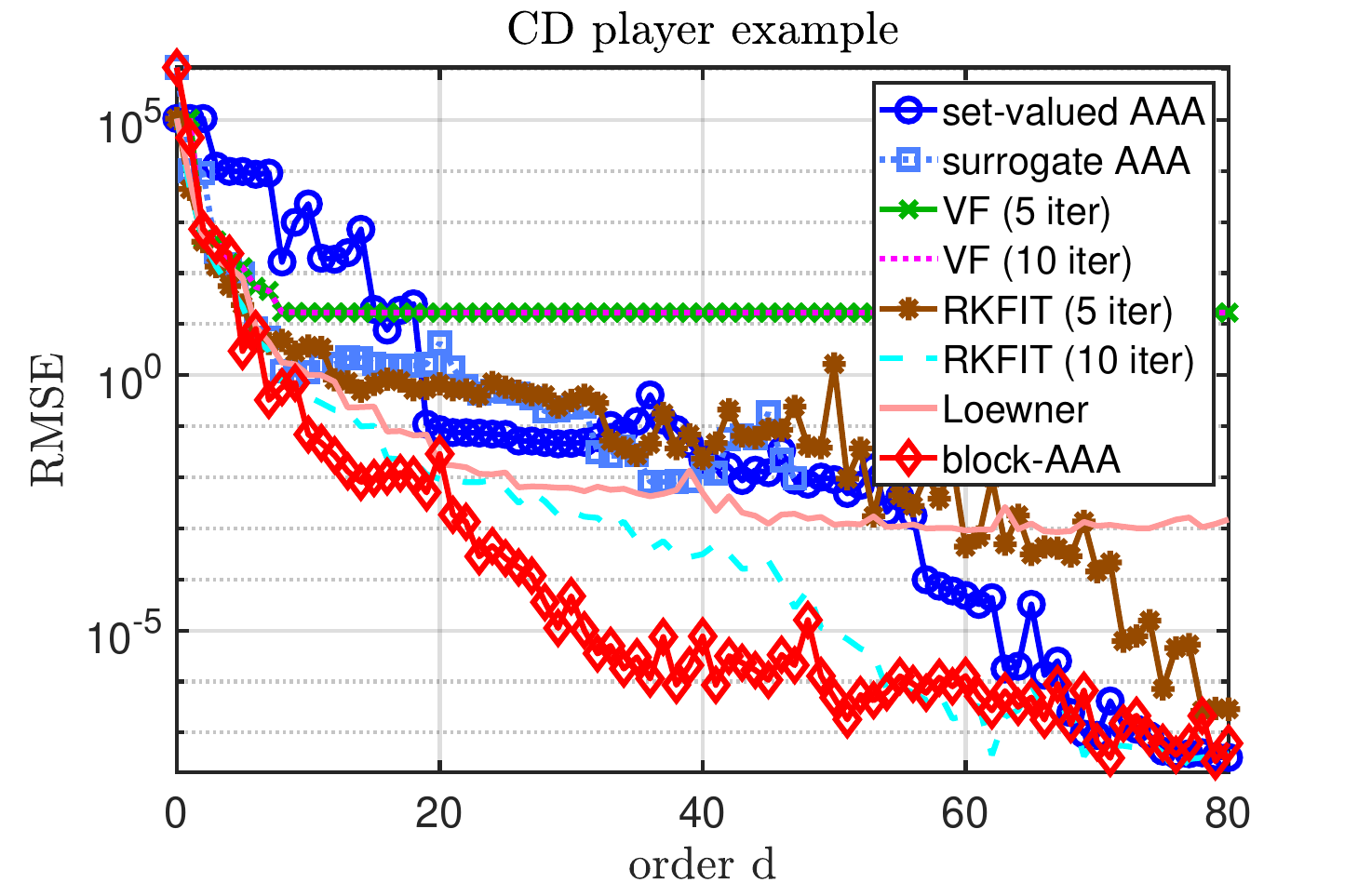}
		\vspace{-6mm}
		\caption{\texttt{CD player} example from Section~\ref{subsec:CD}. Left: entries of the $2\times 2$  matrix $F(z)$ evaluated at 200 points on the imaginary axis. Right: accuracy performance.}
		\label{fig:example_mftex33}
	\end{figure}

	The system's transfer function is a rational function of size $2 \times 2$ that is sampled at $200$ logarithmically-spaced points in the interval $[10^1,10^5]\mathrm{i}$. Figure~\ref{fig:example_mftex33} (left) shows the magnitude of the 4 matrix entries over the frequency range. Figure~\ref{fig:example_mftex33} (right) shows the RMSE values achieved by each of the algorithms for varying orders~$d=0,1,\ldots,80$. The numerical RMSE values and the corresponding timings for orders $d=10$ and $d=20$ are listed in Table~\ref{table_example3}.
	
	\begin{table}[h] 
		\caption{Selected RMSE values and timings for all tested algorithms --- \texttt{CD player} example}
		\label{table_example3}
		\begin{center}
			\begin{tabular}{|l|c|c|c|c|c|l|}\hline
				& \multicolumn{2}{c|}{RMSE} & \multicolumn{2}{c|}{Runtime (ms)}\\
				& $d=10$ & $d=20$ & $d=10$ & $d=20$  \\\hline
				Set-valued AAA & $2.258 \cdot 10^{3}$   & $8.564 \cdot 10^{-2}$ &  $9.0$   & $14.4$  \\
				Surrogate AAA  & $1.129 \cdot 10^{0}$    & $4.301 \cdot 10^{0}$ &  $5.5$   & $16.7$ \\
				VF (5 iter) &  $1.702 \cdot 10^{1}$ &  $1.675 \cdot 10^{1}$ &   $16.9$   & $33.6$ \\
				VF (10 iter) &  $1.683 \cdot 10^{1}$   & $1.674 \cdot 10^{1}$ &  $22.5$   & $37.0$  \\
				RKFIT (5 iter) &  $3.737 \cdot 10^{0}$   & $5.270 \cdot 10^{-1}$ &  $42.4$   & $93.4$  \\
				RKFIT (10 iter) &  $3.806 \cdot 10^{-1}$   & $9.061 \cdot 10^{-3}$ &  $68.0$   & $141.9$  \\
				Loewner &  $1.011 \cdot 10^{0}$   & $1.784 \cdot 10^{-2}$ &  $9.4$   & $10.6$ \\
				block-AAA &  $6.897 \cdot 10^{-2}$   & $2.863 \cdot 10^{-2}$  &  $43.9$   & $96.6$ \\
				\hline
			\end{tabular}
		\end{center}
	\end{table}
	
	\subsection{A second model order reduction example} \label{subsec:ISS}
\begin{figure}[h]
		\includegraphics[scale=.33]{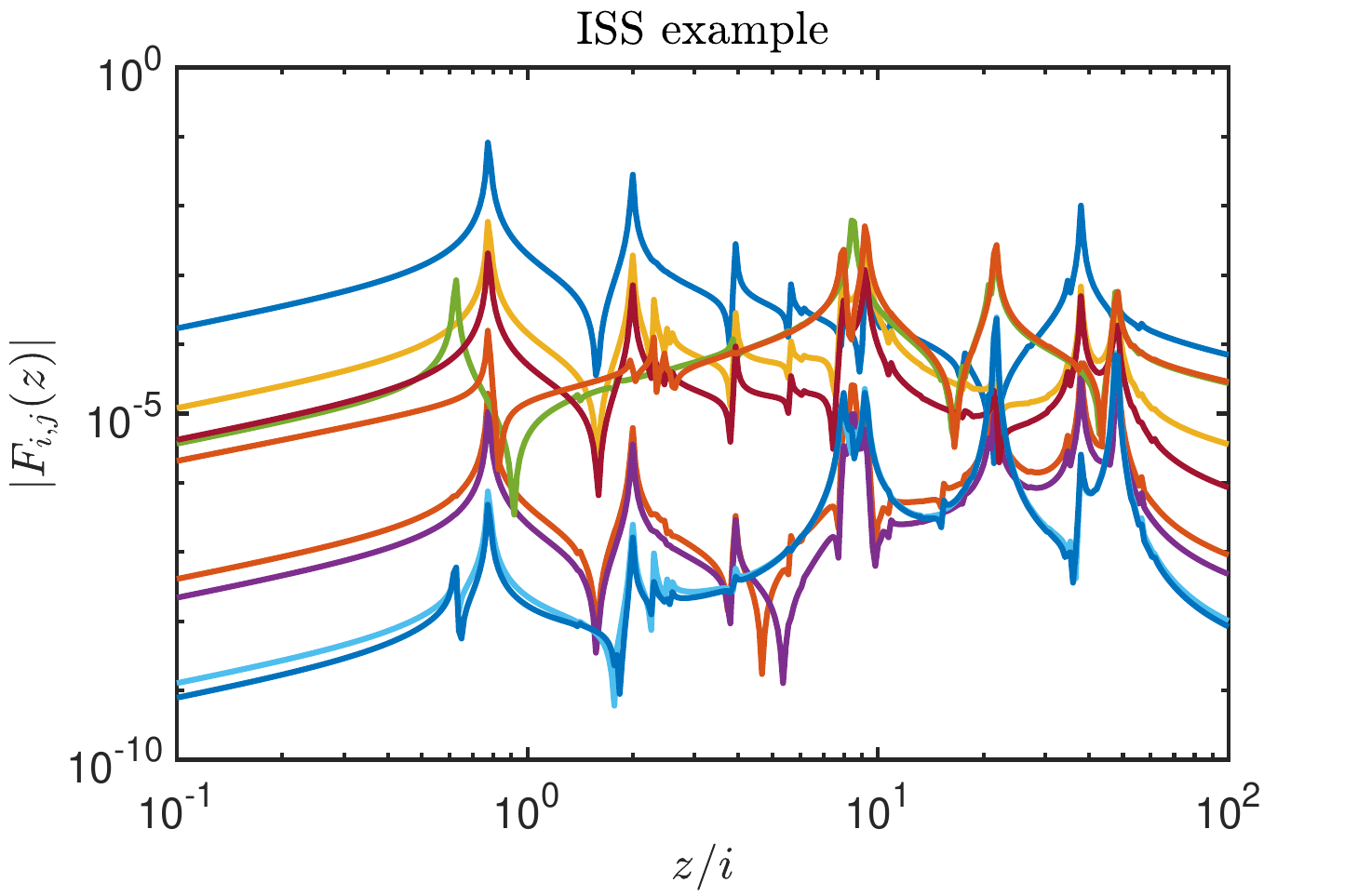}\hspace*{-3mm}\includegraphics[scale=.33]{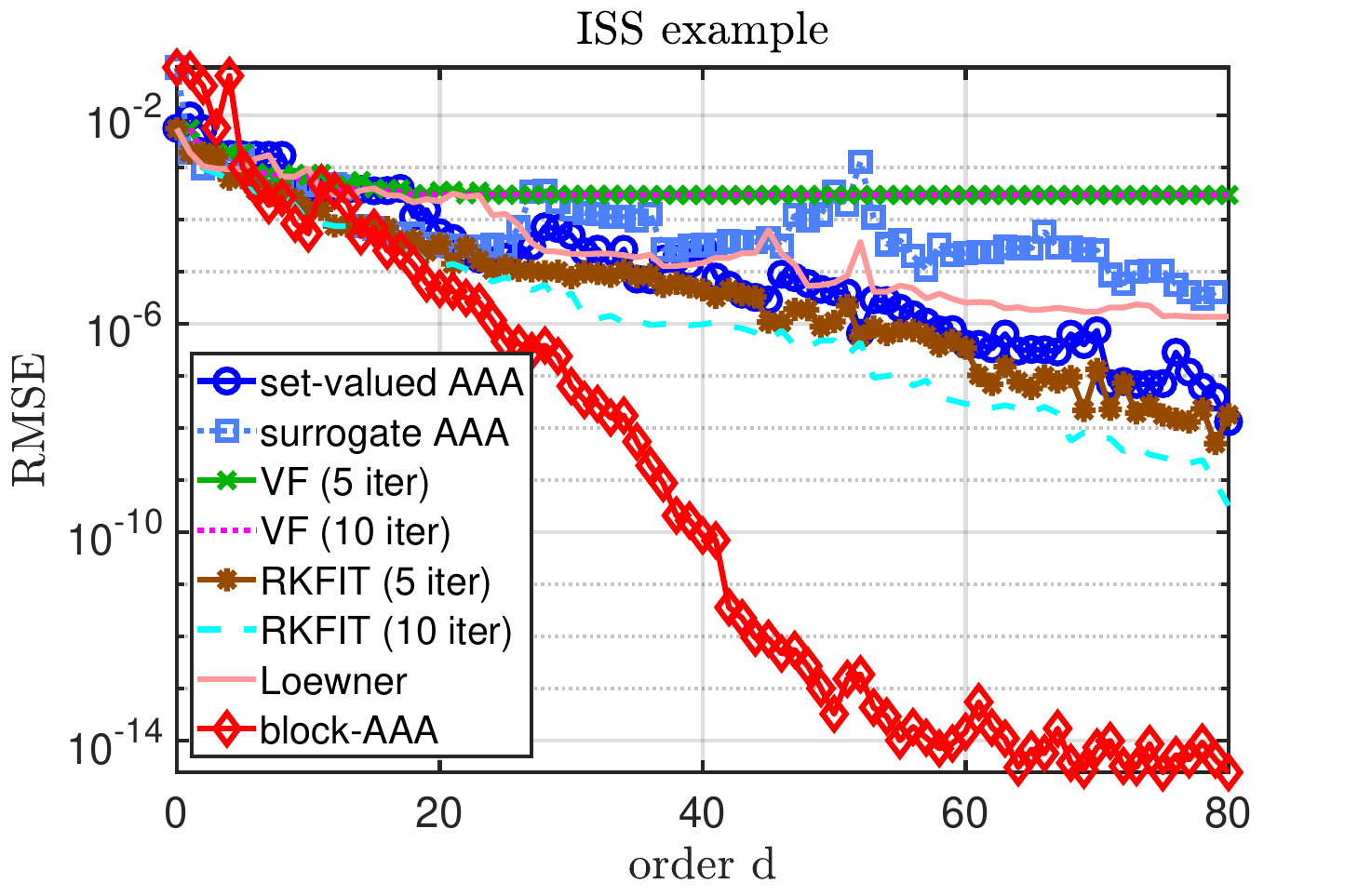}
		\vspace{-6mm}
		\caption{\texttt{ISS} example from Section~\ref{subsec:ISS}. Left: entries of the $3\times 3$ matrix $F(z)$ evaluated at 400 points on the imaginary axis. Right: accuracy performance.}
		\label{fig:example_mftex44}
		\vspace{-0mm}
	\end{figure}

	We consider another model reduction example, namely the ISS model from the SLICOT benchmark collection~\cite{bench02}. 
	There, an LTI system is used as structural  model for the component 1R (Russian service module) in the International Space Station (ISS). The state space dimension of the linear system is $270$ with $m = n= 3$ inputs and outputs. The matrix transfer function of this system is sampled at 400 logarithmically spaced points  in the interval $[10^{-1},10^2]\mathrm{i}$. In Figure~\ref{fig:example_mftex44} (left) we depict the absolute value of the matrix entries. Figure~\ref{fig:example_mftex44} (right) shows the RMSE achieved by each of the algorithms for varying orders $d=0,1,\ldots,80$, and selected numerical RMSE values and corresponding timings for orders $d=10$ and  $d=20$ are given in Table~\ref{table_example4}.
	
	\begin{table}[h] 
		\caption{Selected RMSE values and timings for all tested algorithms --- \texttt{ISS} example}
		\label{table_example4}
		\begin{center}
			\begin{tabular}{|l|c|c|c|c|c|l|}\hline
				& \multicolumn{2}{c|}{RMSE} & \multicolumn{2}{c|}{Runtime (ms)}\\
				& $d=10$ & $d=20$ & $d=10$ & $d=20$  \\\hline
				Set-valued AAA & $3.895 \cdot 10^{-4}$   & $5.543 \cdot 10^{-5}$ &  $12.7$   & $25.3$  \\
				Surrogate AAA  & $4.662 \cdot 10^{-4}$    & $3.209 \cdot 10^{-5}$ &  $8.8$   & $15.7$ \\
				VF (5 iter) &  $6.811 \cdot 10^{-4}$ &  $3.523 \cdot 10^{-4}$ &   $26.7$   & $48.7$ \\
				VF (10 iter) &  $6.729 \cdot 10^{-4}$   & $3.016 \cdot 10^{-4}$ &  $42.6$   & $95.6$  \\
				RKFIT (5 iter) &  $1.555 \cdot 10^{-4}$   & $3.472 \cdot 10^{-5}$ &  $74.7$  & $165.3$  \\
				RKFIT (10 iter) &  $8.735 \cdot 10^{-5}$   & $1.253 \cdot 10^{-5}$ &  $109.4$   & $252.3$  \\
				Loewner &  $9.419 \cdot 10^{-4}$   & $2.225 \cdot 10^{-4}$ &  $23.0$   & $31.2$ \\
				block-AAA &  $5.378 \cdot 10^{-5}$   & $4.678 \cdot 10^{-6}$  &  $102.2$   & $251.7$ \\
				\hline
			\end{tabular}
		\end{center}
	\end{table}
	
	\subsection{A nonlinear eigenvalue problem} \label{subsec:buck}
\begin{figure}[h]
		\includegraphics[scale=.33]{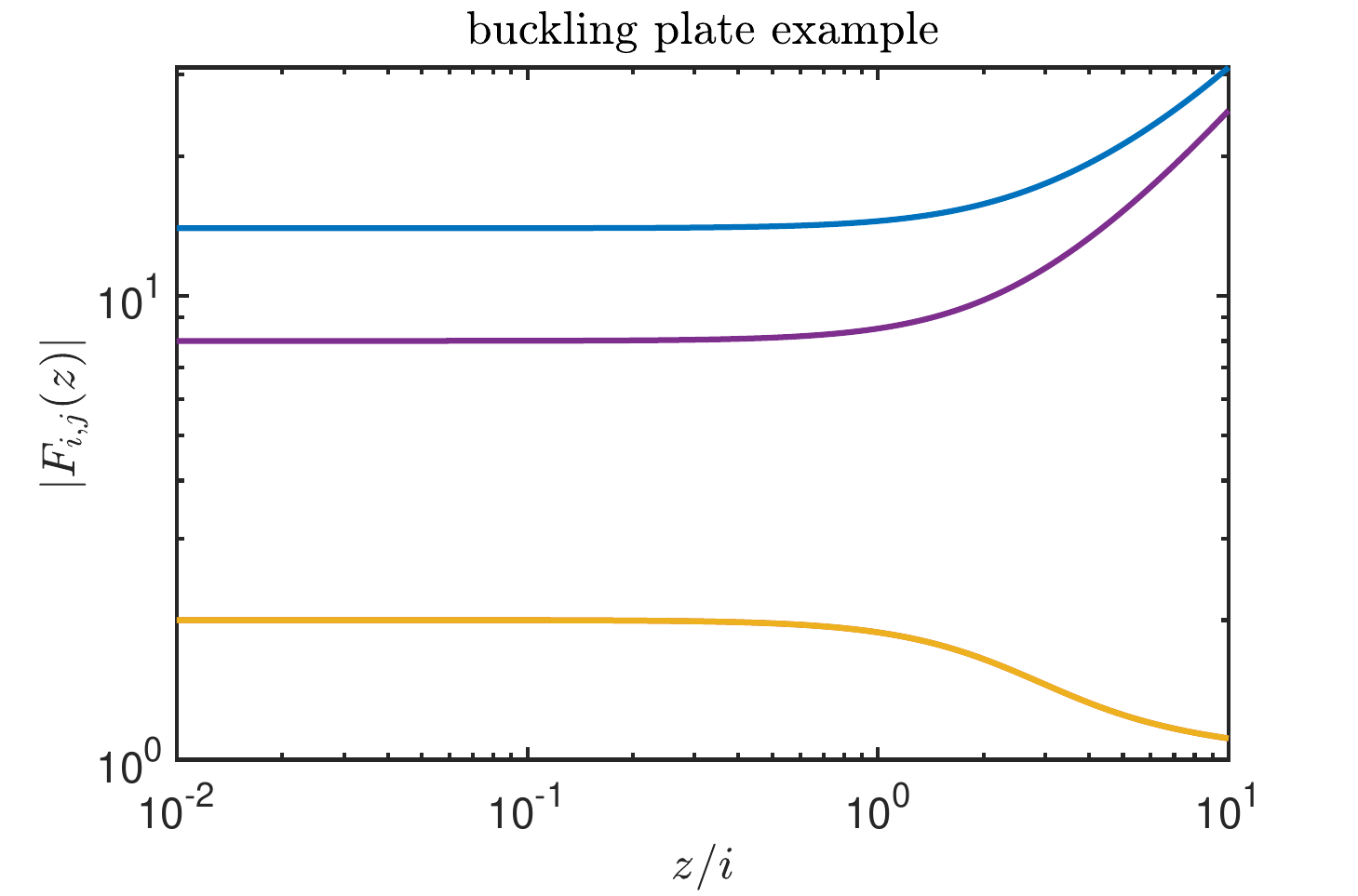}\hspace*{-3mm}
		\includegraphics[scale=.33]{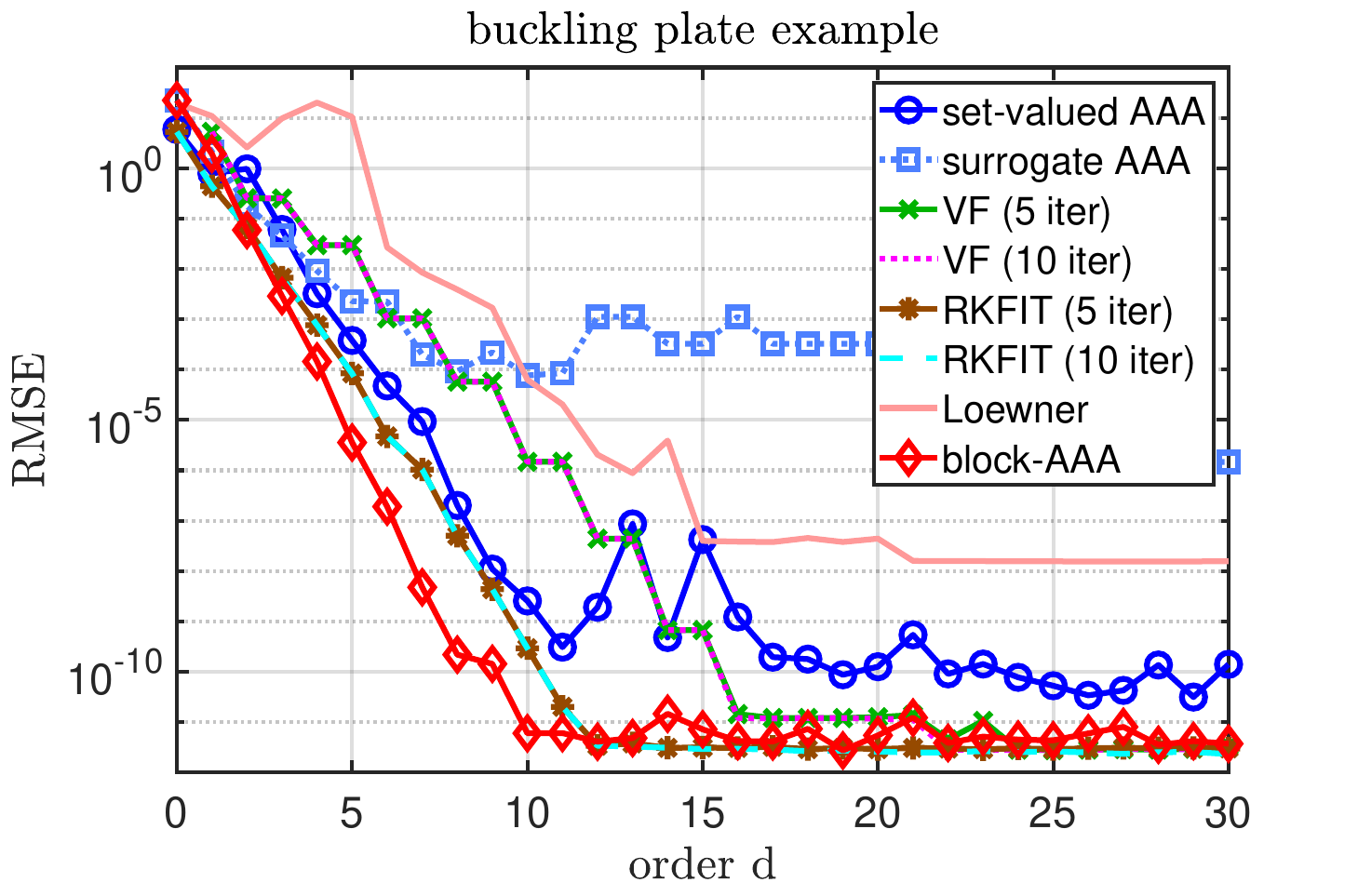}
		\vspace{-6mm}
		\caption{The \texttt{buckling}  problem in Section~\ref{subsec:buck}. Left: entries of the $2\times 2$ (symmetric) matrix $F(z)$ evaluated at 500 points on the imaginary axis. Right: accuracy performance.}
		\label{fig:example_mftex5}
		\vspace{-0mm}
	\end{figure}

	We consider the \texttt{buckling}  example in \cite{HNT19}, a $3 \times 3$ nonlinear eigenvalue problem that arises from a buckling plate model. Since we are interested in the approximation of non-rational functions by means of rational functions, we select only the non-constant part of $F(z)$. Hence, we consider the following $2 \times 2$ symmetric matrix-valued function 
	\begin{equation*}
	F(z) = \left[ \begin{matrix}
	\frac{z(1-2z\cot{2z})}{\tan(z)-z} +10 &  \frac{z(2\lambda-\sin{2z})}{\sin(2z)(\tan(z)-z)} \\[2mm]
	\frac{z(2\lambda-\sin{2z})}{\sin(2z)(\tan(z)-z)} &  \frac{z(1-2z\cot{2z})}{\tan(z)-z} +4
	\end{matrix} \right].
	\end{equation*}
	We choose 500 logarithmically-spaced sampling points  in the interval $[10^{-2},10]\mathrm{i}$. Figure~\ref{fig:example_mftex5} shows the RMSE values achieved by the tested algorithms for varying orders $d=0,1,\ldots,30$. As before, we report the numerical RMSE values and the corresponding timings for two selected orders  $d=10$ and $d=20$ in Table~\ref{table_example5}.
	
	\begin{table}[h] 
		\caption{Selected RMSE values and timings for all tested algorithms --- \texttt{buckling} problem}
		\label{table_example5}
		\begin{center}
			\begin{tabular}{|l|c|c|c|c|c|l|}\hline
				& \multicolumn{2}{c|}{RMSE} & \multicolumn{2}{c|}{Runtime (ms)}\\
				& $d=10$ & $d=20$ & $d=10$ & $d=20$  \\\hline
				Set-valued AAA & $2.543 \cdot 10^{-9}$   & $1.257 \cdot 10^{-10}$ &  $14.9$   & $22.3$  \\
				Surrogate AAA  & $7.640 \cdot 10^{-5}$    & $3.259 \cdot 10^{-4}$ &  $21.6$   & $37.5$ \\
				VF (5 iter) &  $1.476 \cdot 10^{-6}$ &  $1.248 \cdot 10^{-11}$ &   $31.3$   & $46.4$ \\
				VF (10 iter) &  $1.476 \cdot 10^{-6}$   & $1.162 \cdot 10^{-11}$ &  $36.4$   & $68.3$  \\
				RKFIT (5 iter) &  $2.924 \cdot 10^{-10}$   & $2.905 \cdot 10^{-12}$ &  $71.7$   & $160.3$  \\
				RKFIT (10 iter) &  $2.924 \cdot 10^{-10}$   & $2.599 \cdot 10^{-12}$ &  $105.8$  & $ 241.7$  \\
				Loewner &  $6.309 \cdot 10^{-5}$   & $4.397 \cdot 10^{-8}$ &  $31.2$   & $42.3$ \\
				block-AAA &  $5.988 \cdot 10^{-12}$   & $5.438 \cdot 10^{-12}$  &  $104.4$   & $239.9$ \\
				\hline
			\end{tabular}
		\end{center}
	\end{table}
	
	\subsection{A scalar example with noise} \label{subsec:noise0}
	
	In this experiment we investigate the effects of noisy perturbations on the approximation quality. As a test case, we use a scalar function $f(z) = (z-1)/(z^2+z+2)$. \rev{We consider the  RKFIT and  VF methods (non-interpolatory) in comparison to the AAA and Loewner methods (interpolatory).} 
	We sample the function $f(z)$ at 500 logarithmically spaced points in the interval $[10^{-1},10^1]\mathrm{i}$, and then add normally distributed noise with a  standard deviation of $\tau = 10^{-2}$ to these samples. 
	
\begin{figure}[h]
		\includegraphics[scale=.33]{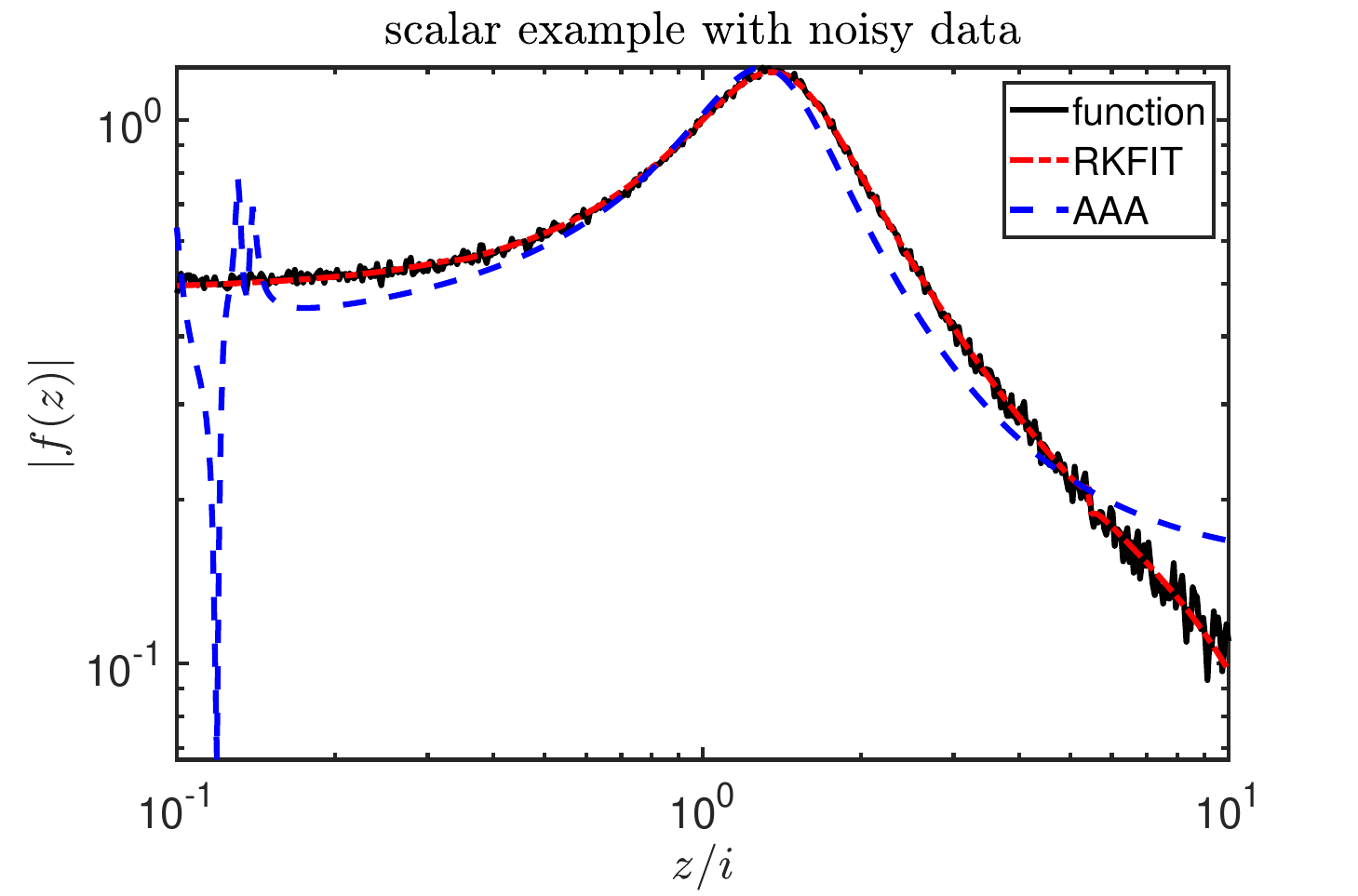}\hspace*{-3mm}
		\includegraphics[scale=.33]{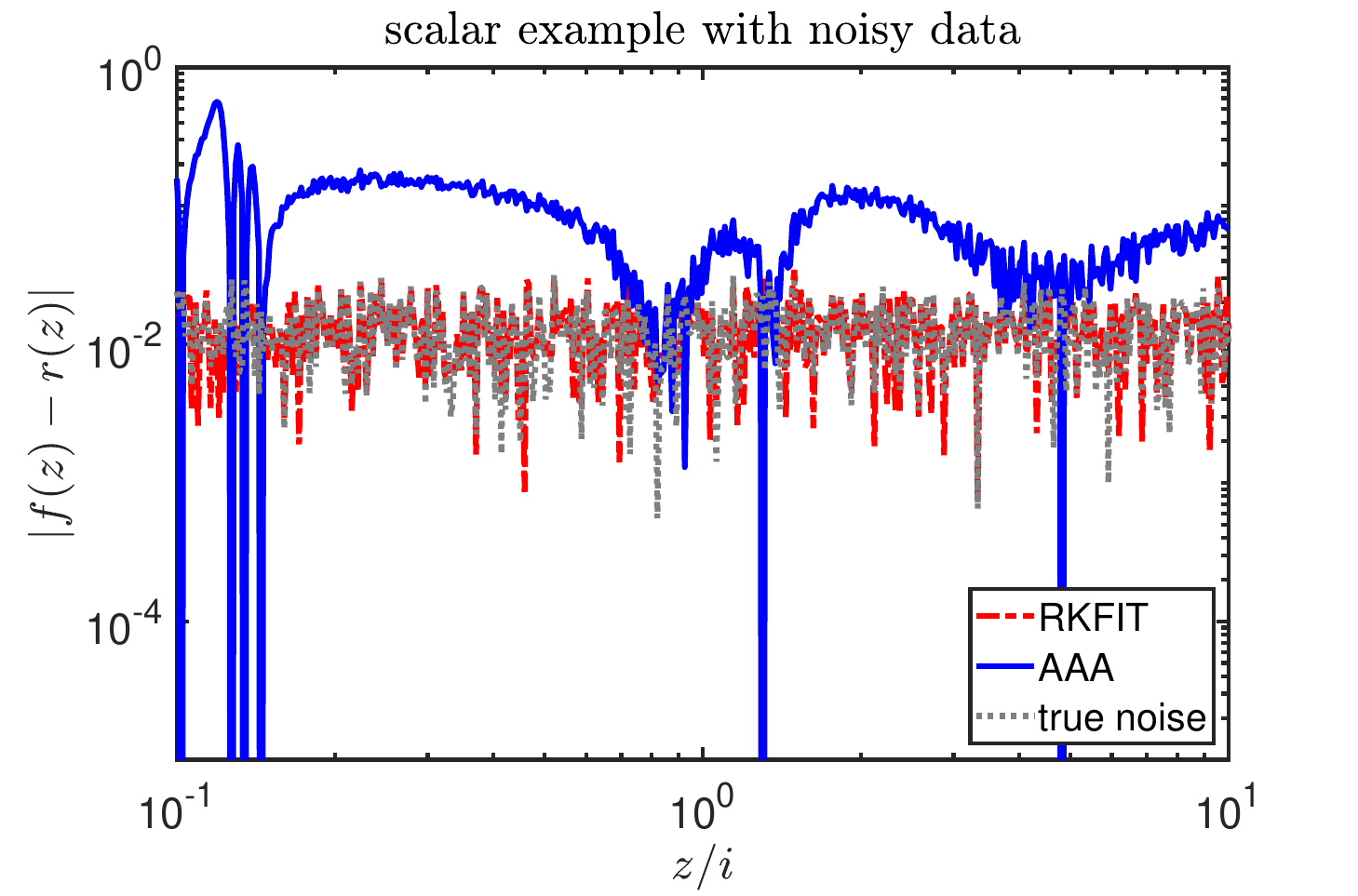}
		\vspace{-6mm}
		\caption{Scalar example with noisy measurements from Section~\ref{subsec:noise0}. Left: the \rev{perturbed version of the} scalar function $f(z)$ and the computed approximants of order $d=5$. Right: the deviation between the RKFIT and AAA approximants and the original function $f(z)$.}
		\label{fig:example_mftex77}
		\vspace{-0mm}
	\end{figure}

	We first compute rational approximants of degree $d=5$ using the above methods. For RKFIT, we perform three iterations starting with the default initialization of all poles at infinity. As it can be observed in Figure~\ref{fig:example_mftex77} (left), the red curve corresponding to the RKFIT approximants follows the noisy measurements very well on average. At the same time, the blue curve corresponding to the AAA approximant shows considerable deviations from the measurements.  By inspecting the deviation between the two rational approximants and the original function $f(z)$ in Figure~\ref{fig:example_mftex77} (right), we find that the RKFIT approximation error is comparable to the noise level. The RKFIT approximant has effectively estimated the additive noise rather accurately, which is also confirmed visually by the true noise curve overlaid as a dotted line on top of the RKFIT error curve. \rev{For this scalar problem, the VF approximant behaves very similarly to the RKFIT approximant, achieving a stagnation RMSE close to the noise level.} On the other hand, the approximation error attained by the AAA method is between 1--2 orders of magnitude larger than the added noise level. \rev{The Loewner method stagnates on an even higher level.} 
	
\begin{figure}[h]
		\centering\includegraphics[scale=.38]{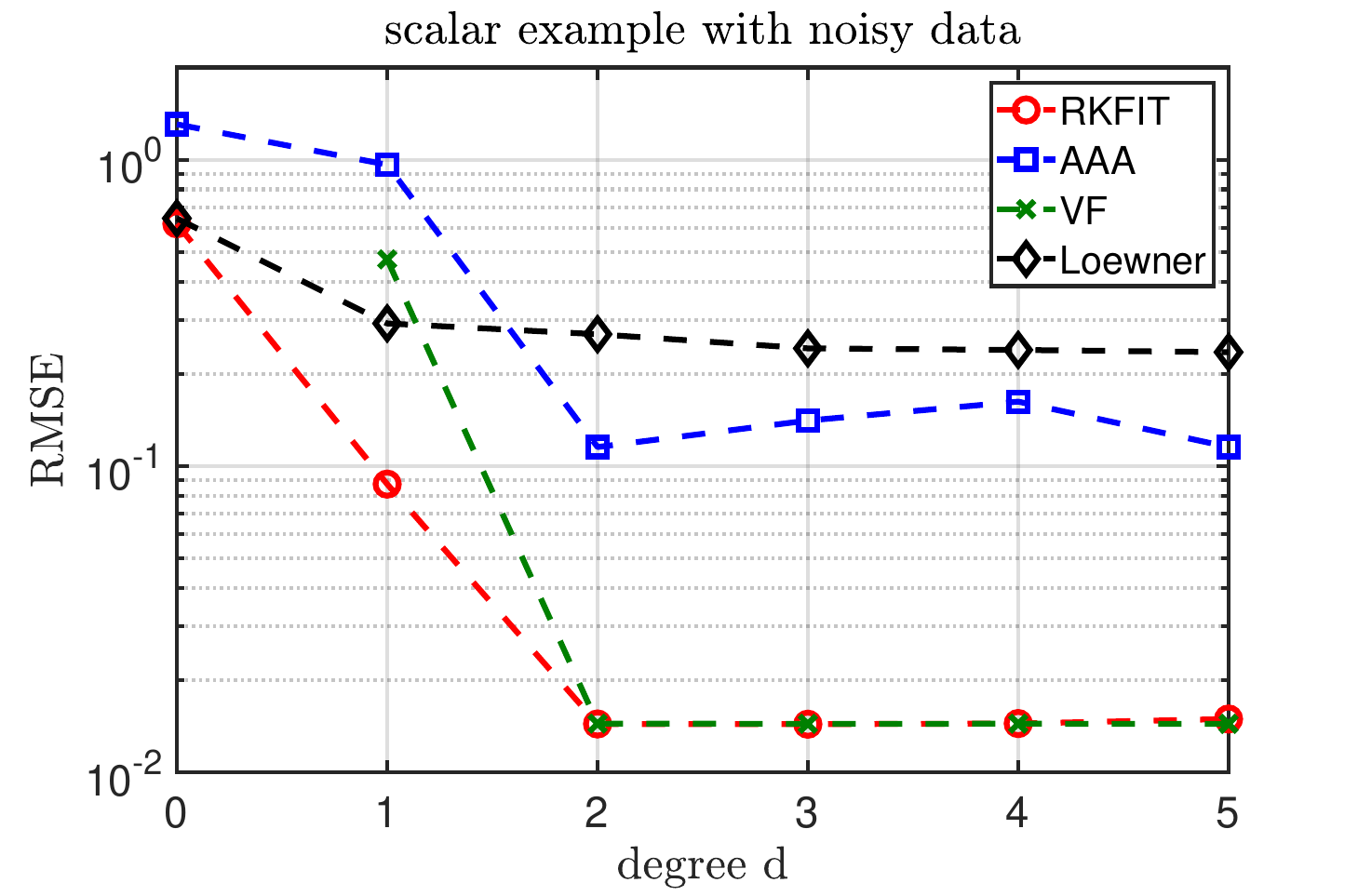}
		\vspace{-2mm}
		\caption{\rev{RMSE convergence of the approximants computed via non-interpolatory (RKFIT, VF) and  interpolatory (AAA, Loewner) methods for increasing degree $d$.} The original function is rational of degree~2, and noise with a standard deviation of $\tau=10^{-2}$ has been added to perturb it.}
		\label{fig:example_mftex7}
		\vspace{0mm}
	\end{figure}

	Finally, we vary the degree $d=0,1,\ldots,5$ and compute the corresponding RMSE values of the RKFIT and AAA approximants. The results are depicted in Figure~\ref{fig:example_mftex7}. Note that the RKFIT method achieves lower RMSE values  than the AAA method and exhibits a more "regular" convergence. The RKFIT error stagnates at the noise level of approximately $\tau = 10^{-2}$ when the degree $d=2$ is reached. We believe that this preferable approximation performance is due to the non-interpolatory nature of RKFIT, and that AAA suffers from the fact that noisy values are being interpolated (resulting in the observed oscillations of the approximants). 
	
	\subsection{A second example with noise} \label{subsec:noise}
	
	For this case, we analyze a noisy version of the \texttt{ISS} model considered in Section~\ref{subsec:ISS}. We modify the original LTI system $(A,B,C)$ by choosing $C = B^T$ in order to obtain a symmetric transfer function and be able to apply the Matrix Fitting Toolbox. 
	
\begin{figure}[h]
		\includegraphics[scale=.33]{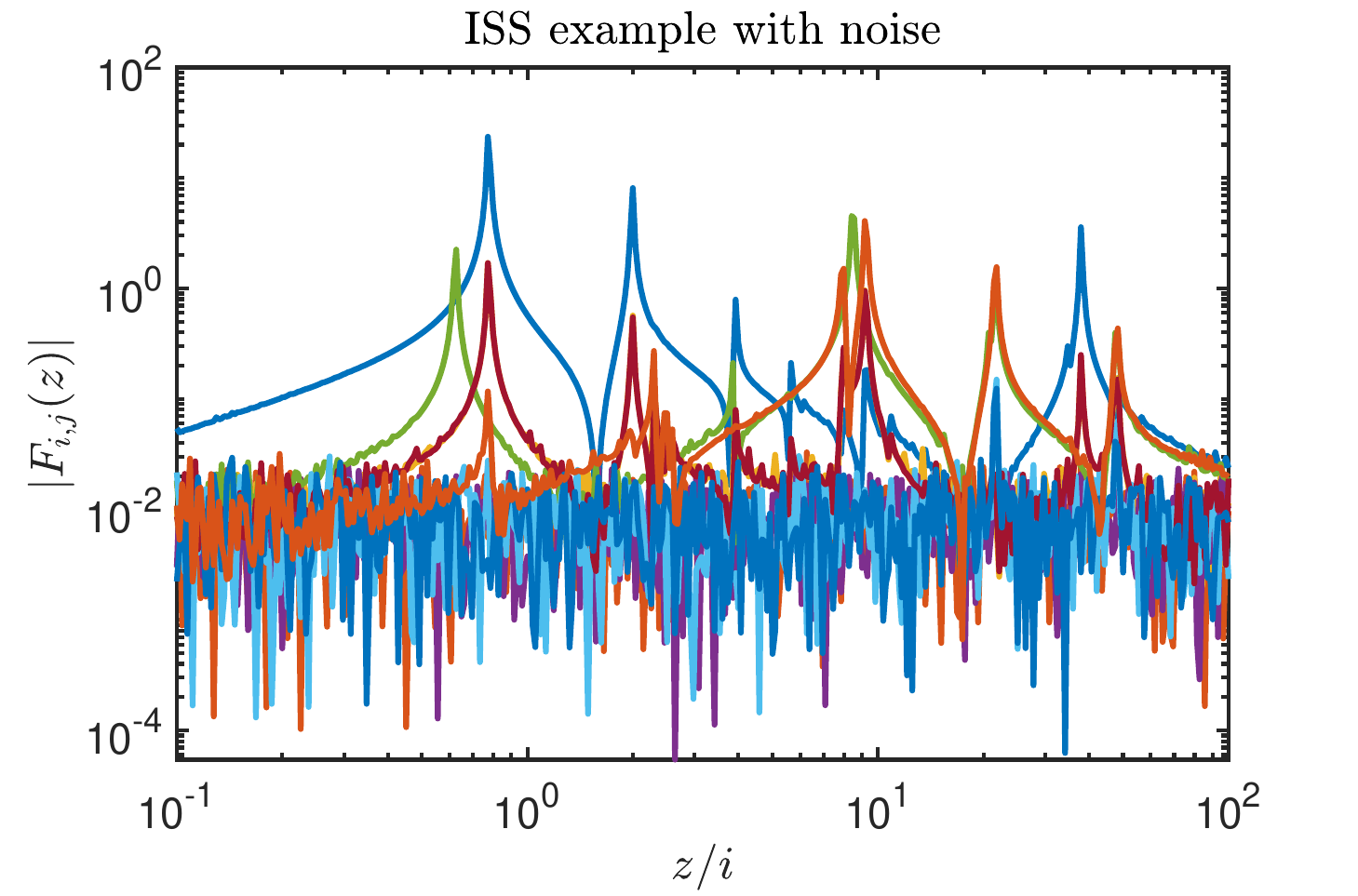}
		\hspace*{-3mm}\includegraphics[scale=.33]{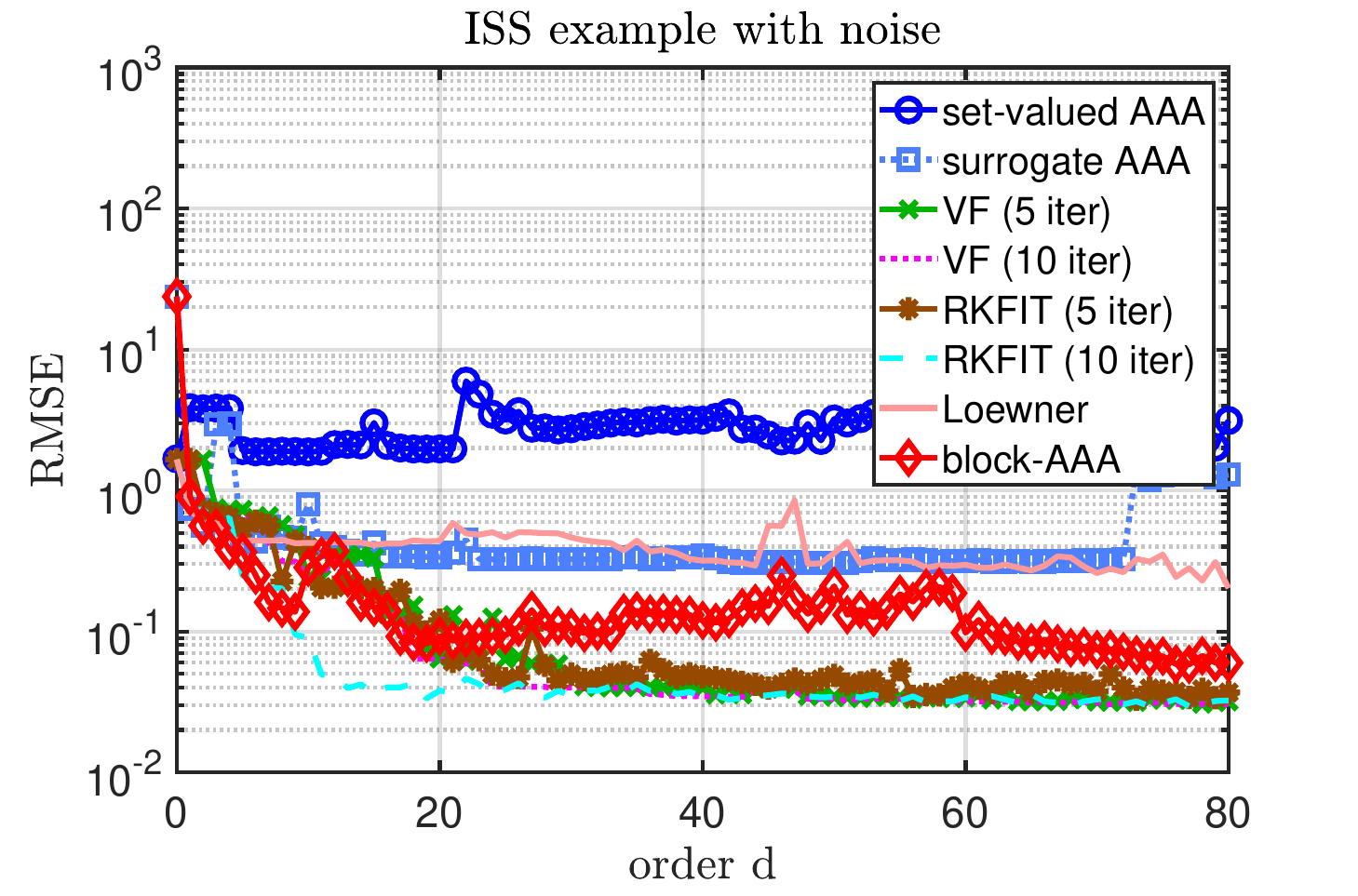}
		\vspace{-6mm}
		\caption{Noisy \texttt{ISS} example from Section~\ref{subsec:noise}. Left: the  entries of the perturbed $3\times 3$ matrix~$F(z)$. Right: performance comparison.}
		\label{fig:example_mftex66}
		\vspace{0mm}
	\end{figure}

	The 400 sampling values are corrupted by additive normally distributed noise with standard deviation $\tau = 10^{-2}$. In Figure~\ref{fig:example_mftex66} (left) we depict the magnitude of each entry of the perturbed $3 \times 3$ matrix transfer function. 
	
	The RMSE values for varying orders $d=0,1,\ldots,80$ are plotted in Figure~\ref{fig:example_mftex66} (right). 
	Now, for $d=10$ and $d=20$ we record the numerical RMSE values for all methods in Table~\ref{table_example6}. In accordance with our observations on the previous scalar example (Section~\ref{subsec:noise0}) we find that the non-interpolatory methods like RKFIT and VF exhibit the most steady convergence behavior, while the interpolation-based methods (AAA, Loewner) produce approximants whose accuracy varies wildly as the order $d$ changes. The only methods that reliably attain an RMSE close to the noise level of $\tau=10^{-2}$ are RKFIT and VF, with RKFIT (10 iterations) consistently attaining the lowest RMSE values for all considered orders.
	
	\rev{Recall that the original (nonsymmetric) \texttt{ISS} model was treated (without adding noise) in Section~\ref{subsec:ISS}. As it can be seen in Table~\ref{table_example6}, the timings of all methods are similar to those reported in Table~\ref{table_example4} for the case without noise.}
	
	\begin{table}[h] 
		\caption{Selected RMSE values and timings for all tested algorithms - noisy \texttt{ISS} example}
		\vspace{-0mm}
		\label{table_example6}
		\begin{center}
			\begin{tabular}{|l|c|c|c|c|c|l|}\hline
				& \multicolumn{2}{c|}{RMSE} & \multicolumn{2}{c|}{Runtime (ms)}\\
				& $d=10$ & $d=20$ & $d=10$ & $d=20$  \\\hline
				Set-valued AAA &  $1.877 \cdot 10^{0}$   & $1.971 \cdot 10^{0}$ & 15.2 & 29.3 \\
				\hline 
				Surrogate AAA  & $7.994 \cdot 10^{-1}$    & $3.380 \cdot 10^{-1}$ & 11.4 & 23.3 \\
				\hline
				VF(5 iter) &  $2.388 \cdot 10^{-1}$ &  $6.417 \cdot 10^{-2}$ & 33.2 & 56.8 \\
				\hline
				VF (10 iter) &  $2.387 \cdot 10^{-1}$   & $6.404 \cdot 10^{-2}$ & 46.9 & 91.4 \\
				\hline
				RKFIT (5 iter) &  $3.946 \cdot 10^{-1}$   & $1.201 \cdot 10^{-1}$ & 56.7 & 126.6 \\
				\hline
				RKFIT (10 iter) &  $8.991 \cdot 10^{-2}$   & $3.810 \cdot 10^{-2}$ & 109.1 & 262.3 \\
				\hline
				Loewner &  $4.250 \cdot 10^{-1}$   & $4.396 \cdot 10^{-1}$ & 21.5 & 30.0 \\
				\hline
				block-AAA &  $2.847 \cdot 10^{-1}$   & $9.215 \cdot 10^{-2}$ & 102.4 & 255.4 \\
				\hline
			\end{tabular}
		\end{center}
	\end{table}

	\section{Discussion}
	\label{sec:discussion}
	
	
	In this section we present a detailed analysis of the numerical results reported in Section~\ref{sec:numer}.
	The discussion will address each method separately, taking into consideration the approximation quality, the amount of time needed to run the method (runtime), and the time needed to evaluate the computed approximants (evaltime).
	
	\subsection{Set-valued AAA}
	An advantage of the set-valued AAA algorithm is its fast runtime. With some exceptions, the set-valued AAA method was the second fastest method in our tests, surpassed only by the surrogate AAA method (and sometimes by the Loewner method).
	
	With respect to the approximation quality, this method produced models with similar accuracy to RKFIT for the example in Section~\ref{subsec:matfit} (see Figure~\ref{fig:example_mftex22}, right). A similar behavior was observed for the next three examples, producing good and very good results. For the example in Section~\ref{subsec:CD} it even produced results comparable in accuracy to the block-AAA method when $d > 65$ (see Figure~ \ref{fig:example_mftex33}, right). For the examples in Sections~\ref{subsec:noise0} and \ref{subsec:noise} with added noise, however, the method produced poor results. In particular for the noisy ISS example, the RMSE values produced by the set-valued AAA method are the highest (see Figure~\ref{fig:example_mftex66}, right). As previously discussed, the reason for this poor approximation is the interpolation of noisy data values. \rev{Moreover, the evaluation time for the computed approximants was found to be small, comparable to the evaluation of surrogate AAA (the lowest) approximants and those obtained with the Loewner method (see Table~\ref{table_example2}). The orders needed to achieve selected RMSE values were relatively low, only slightly higher than the order of the corresponding RKFIT approximants (see Table~\ref{table_example22}). Finally, the method required the smallest runtimes (by far) for low target accuracy $\epsilon_3 = 10^{-9}$ and $\epsilon_4 = 10^{-12}$; see Table~\ref{table_example222}}.

	\subsection{Surrogate AAA}
	The surrogate AAA method performs best in terms of runtime, and it could be an attractive approach when the problem dimension $(m,n)$ is very large and approximation accuracy is not the main concern. The algorithm produces, in general, quite poor approximants compared to the other methods.  
	\rev{For the example in Section~\ref{subsec:matfit}, surrogate AAA did not reach an RMSE below  $\epsilon_3 = 10^{-9}$, while all other methods were able to reach an RMSE of $\epsilon_4 = 10^{-12}$; see Figure~\ref{fig:example_mftex22} (right) and Table~\ref{table_example22}. In terms of runtime, the surrogate AAA was indeed the fastest for approximants with  RMSE values $\epsilon_1 = 10^{-3}$ and $\epsilon_2 = 10^{-6}$ (see Table~\ref{table_example222})}. A poor approximation was observed for both MOR examples in Section~\ref{subsec:CD} and \ref{subsec:ISS}. For the latter, the method was the second worst (after VF). For the buckling plate example in Section~\ref{subsec:buck}, surrogate AAA produced the poorest approximation results (see Figure\;\ref{fig:example_mftex5}). Also, for the ISS example with noise in Section~\ref{subsec:noise}, we observed that the method produced poor results, comparable only to those returned by the Loewner framework.
	
	
	\subsection{Vector Fitting}
	
	For this method, we have used precisely the implementation provided in~\cite{MatVF}. One limitation of this toolbox is that it can only be applied for symmetric matrix-valued  functions, i.e., the samples must be symmetric matrices. Consequently, this implementation of VF could not be applied for two of the examples, namely the ones in Sections~\ref{subsec:CD} and \ref{subsec:ISS}.
	
	For the example in Section~\ref{subsec:matfit} we observed that the VF method produced poor approximants compared to some of the other methods. \rev{As can be seen in Figure~\ref{fig:example_mftex22} (right) and in Table~\ref{table_example22}, VF was the method that needed the largest orders to achieve the four target RMSE values $\epsilon_k$. Additionally, the corresponding runtimes presented in Table~\ref{table_example222} were the second highest, only better than the ones for RKFIT.}
	
	On the other hand, for the buckling plate example in Section~\ref{subsec:buck}, the VF method surpassed three methods and reached an RMSE of about $10^{-12}$ for order $d=30$ (see Figure~\ref{fig:example_mftex5}). Note also that VF together with RKFIT are the only methods that reliably attain an RMSE value ($\eta=0.03$) close to that of the noise level ($\tau=10^{-2}$), as could be observed in Figure~\ref{fig:example_mftex66} (right). \rev{A similar conclusion can be drawn for the experiments from Section~\ref{subsec:noise0}, e.g., displayed in Figure~\ref{fig:example_mftex66}. There, VF and RKFIT performed better; the true degree equal to 2 was indicated as the degree for which the RMSE values stagnate at a value close to the standard deviation of the noise. Both interpolation-based methods, i.e., AAA and Loewner, failed to do that.}

	Another observation is that it is generally not worth performing 10 VF iterations instead of just 5. We did not observe any significant approximation enhancement and the total runtime of VF will approximately double for 10 VF iterations. Throughout the experiments performed in Section~\ref{sec:numer}, we found that VF was generally slower than the set-valued and surrogate AAA methods, but faster than the RKFIT and block-AAA methods.
	
	
	\subsection{RKFIT}
	
	For the RKFIT method, we have used the implementation in the RKToolbox. For each  experiment, we ran the RKFIT algorithm for a fixed prescribed number of iterations.  
	\rev{When stopping after ten iterations, RKFIT was generally the slowest method in all experiments performed (with runtimes comparable to those of block-AAA). When stopping after five iterations, RKFIT was generally the third slowest method. Additionally, it required one of the highest evaluation times (see Table~\ref{table_example2}) and it needed the highest runtimes for achieving each of the target tolerances $\epsilon_k$ in Table~\ref{table_example222}. 
		
		On the other hand, in all cases listed in Table~\ref{table_example22}, RKFIT achieved  RMSE tolerance values with the second-lowest order; only block-AAA required lower orders.} Moreover,  RKFIT  managed to yield the second lowest RMSE values for the first five experiments, surpassed only by block-AAA. \rev{This is not surprising given that RKFIT produces approximants $R_d$ where all matrix entries share a common scalar denominator of degree $d$, say, while a block-AAA approximant of the same order~$d$ can have up to $dm$ singularities. } 
	Hence, among the scalar denominator methods, RKFIT was the  most accurate and reliable. Moreover, RKFIT proved to be the most robust when dealing with noise. For both  variants (with 5 and 10 iterations), the RMSE values obtained by RKFIT were the lowest among all methods (see Figure~\ref{fig:example_mftex66}, right), comparable only to those of VF.
	
	The runtime of the RKFIT method linearly depends on the number of iterations.  It was observed in some cases that the RMSE values would not significantly decrease when increasing the number of RKFIT iterations from five to ten. For example, it can be seen in Figure~\ref{fig:example_mftex22} (right) and also in Figure~\ref{fig:example_mftex5} (right), that the two RMSE curves for RKFIT (5 iterations) and RKFIT (10 iterations) are practically indistinguishable. \rev{It might have been possible to achieve a similar accuracy with fewer than 5 iterations. In practice, a  stopping criterion based on the stagnation of  RMSE should therefore be used. Indeed, as RKFIT is based on least squares fitting, the RMSE is available at each iteration for ``free'' (computed from a thin SVD which needs to be computed anyway), and hence can be used for a stagnation stopping criterion. Further, it is also possible to have RKFIT determine the degree $d$ dynamically by adding one pole at a time until a target tolerance is reached. This is implemented in RKToolbox when RKFIT is called with an empty list \texttt{xi = []} for the initial poles. We have not used any of these options here as it would have significantly complicated the comparison with the other methods.}
	
	\subsection{The Loewner framework}
	
	The Loewner framework, as introduced in Section \ref{subsec:Loewner_scalar}, is the only ``direct'' method out of the six included in this study, i.e., it does not rely on an iteration. This could be an advantage with respect to the speed of execution. In general, this method performs very well in terms of runtime for small, and medium to large data sets. The reason for this is that it relies on computing a full singular value decomposition of a matrix with dimension half of the data set.
	
	It was observed that, with some exceptions, the Loewner framework was the third fastest method, surpassed only by the set-valued and surrogate AAA methods. This is not surprising since all the three methods rely on multiplying the matrix-valued samples with left and right tangential vectors, thereby reducing the dimension of the SVD problem to be solved at each iteration.
	
	In terms of  approximation accuracy, the Loewner framework produced rather poor results as compared to other methods. For example, in Section~\ref{subsec:matfit} it was observed that the RMSE values achieved by this method were higher than others, only comparable to those of VF (see Figure~\ref{fig:example_mftex22}, right). For the experiments in Sections~\ref{subsec:ISS} and \ref{subsec:buck}, it produced better results only than those of surrogate AAA. Note that, for the CD player example in Section~\ref{subsec:CD}, the approximation quality was comparable to that of the set-valued AAA method and RKFIT (5 iterations). \rev{It is to be noted that the Loewner method produced poor results for the scalar noisy example (see Figure\;\ref{fig:example_mftex7}), and failed for the ISS example with added noise (see Figure\;\ref{fig:example_mftex66}). In terms of evaluation time, the Loewner framework performed well having second and third lowest times (see Table~\ref{table_example2}). Similar observation were made with respect to the runtimes needed to reach a specific target accuracy, presented in Table~\ref{table_example222}---note that only set-valued AAA was faster. Finally, the Loewner method required rather large orders to reach the target RMSEs in Table~\ref{table_example22}.}
	
	\subsection{Block-AAA}
	
	The block-AAA method introduced in Section~\ref{sec:alg} is the only method in this comparison which produces rational approximants with \emph{nonscalar} denominators. As such, the approximation quality obtained for a certain order~$d$ might not be fully comparable to that of another method. In terms of the execution times, the method proved to be the second slowest for the examples in Sections~\ref{subsec:matfit}, \ref{subsec:CD}, and \ref{subsec:buck} (in these cases, only the RKFIT method with 10 iterations was slower). Additionally, for the example presented in Section~\ref{subsec:ISS}, the block-AAA method was indeed the slowest. The reason for this inefficiency is the repeated solution of the SVD problem in Step~5 of the algorithm (see Figure~\ref{alg:baaa}), involving a matrix of size $m(j+1)\times \ell n$ at each iteration $j=0,1,\ldots,d$, as well as the greedy search in Step~2 which requires the evaluation of $R_{j-1}(z)$ at  all remaining sampling points in $\Lambda^{(j)}$ to determine the next interpolation point~$z_j$. We have tried to speed up Step~5 using the updating trick of the SVD decomposition described in~\cite{LPVM18}, but noticed some numerical instabilities for sufficiently large degrees. We have therefore not used any SVD updating strategy for the block-AAA results reported here.
	
	In terms of approximation accuracy, the block-AAA method usually produced the best results for a given order~$d$, with the noisy ISS MOR example being the only exception. \rev{A representative illustration is Figure~\ref{fig:example_mftex22} (right) and also Table~\ref{table_example22}; we note that block-AAA requires order of $d=21$ to reach an RMSE value of $\epsilon_4 = 10^{-12}$, while RKFIT and the set-valued AAA method require an order at least twice higher (other methods require even higher orders or do not reach this value at all). On the other hand, block-AAA approximants are generally very slow to evaluate (see Table~\ref{table_example2}). This is due to their nonscalar denominators, which requires matrix inversions. Additionally, the runtimes needed to reach tolerance values $\epsilon_k$ proved to be lower only than those of RKFIT and VF (10 iterations) and similar to those of VF (5 iterations); see Table~\ref{table_example222}}. A similar behavior was observed in Figure~\ref{fig:example_mftex33} (right) for the CD player example and in Figure~\ref{fig:example_mftex44} (right) for the ISS example. For the noisy data experiment in Section~\ref{subsec:noise}, block-AAA was outperformed by RKFIT and, again, we believe this is due to RKFIT's non-interpolatory character.

	\section{Conclusion}
	\label{sec:conclusion}
	
	We proposed an extension of the AAA algorithm to the case of matrix-valued interpolation. The block-AAA method uses a new generalized barycentric representation of the interpolant with matrix-valued weights. Our numerical experiments indicate that this modification allows for the computation of accurate approximants of low order. Having a low order might be particularly interesting in cases where the approximant needs to be linearized for further processing, e.g., in the context of nonlinear eigenvalue computations. In the appendix we show how the generalized barycentric interpolants can be linearized into a corresponding eigenvalue problem of size proportional to the order. \rev{Hence, for a fixed sampling accuracy, a  linearization obtained from a block-AAA approximant can be of considerably smaller dimension than that  obtained from another scalar-denominator method.} The use of block-AAA for this application will be explored in future work. In terms of algorithmic complexity and runtime, the block-AAA method is currently inefficient and can be practically applied only for problems of small dimensions. Further work is required to improve its performance. One approach to deal with this shortcoming could be to replace the full singular value decomposition performed at each step of block-AAA with a CUR decomposition (selecting only a small number of relevant columns and rows in the block Loewner matrix), and to somehow restrict the search for the next interpolation point to a smaller set. 
	
	Our comparisons indicated that there was no method that always performed best with respect to accuracy and runtime. Interpolation-based methods are often cheaper to run, but they may suffer in the presence of noise for which approximation-based methods perform better. It seems fair to say that, despite the exciting recent developments in the area, there is still a lot of work to be done to design robust and fast methods for the rational approximation of matrix-valued functions. 
	
	\section*{Acknowledgements} \rev{We thank the two anonymous referees for providing many insightful comments that have improved our manuscript.} We would like to thank the Manchester Institute for Mathematical Sciences (MIMS) for enabling a visit of the first author to Manchester during which much of this work was done. We would also like to thank Steven Elsworth, Yuji Nakatsukasa, Gian Maria Negri Porzio, and Nick Trefethen for many discussions on the topic of rational approximation.

	\bibliographystyle{plain}
   \bibliography{references_paper}

	\section*{Appendix: Linearization of generalized barycentric forms} 
	While we have been primarily concerned with evaluating the various algorithms in terms of their RMSE approximation performance, in the context of nonlinear eigenvalue problems one also requires the computed approximants to be suitable for  pole and ``root'' computations. More precisely, in nonlinear eigenvalue problems one is interested in the points $\xi\in\mathbb{C}$ where one or more entries of the rational function $R_d(z)$ have a pole, and the nonlinear eigenvalues $\lambda\in\mathbb{C}$ where $R_d(\lambda)$ is a singular matrix. Naturally the latter problem requires $R_d(z)$ to be a square matrix.

	Let us discuss the suitability of the matrix-valued barycentric forms in Section~\ref{sec:represent} for nonlinear eigenvalue computations based on linearization; see, e.g., \cite[Section~6]{GT17}. 
	First note that \eqref{eq:bb1} and \eqref{eq:bb2} are special cases of \eqref{eq:bb3}, respectively: we obtain \eqref{eq:bb1} from \eqref{eq:bb3} by setting $C_k = w_k F(z_k)$ and $D_k = w_k I$, while we obtain \eqref{eq:bb2} from \eqref{eq:bb3} by setting $C_k = W_k F(z_k)$ and $D_k = W_k$. Hence it is sufficient to focus on the most general representation~\eqref{eq:bb3}. 
	
	If $\lambda\in\mathbb{C}$ is a nonlinear eigenvalue of $R_d(z)$, then clearly $R_d(z)^{-1}$ has a pole at $z=\lambda$, and vice versa. It is therefore sufficient to be able to find the nonlinear eigenvalues of~$R_d$ defined in \eqref{eq:bb3} as those points $\lambda\in\mathbb{C}$ for which the ``numerator''
	\[
	N(z) = \frac{\displaystyle\sum_{k=0}^d \frac{w_k [C_k/w_k]}{z-z_k}}{\displaystyle\sum_{k=0}^\ell \frac{w_k}{z-z_k}}
	\]
	is a singular matrix. It is well known that for given support points $z_0,z_1,\ldots,z_d$, we can choose nonzero weights $w_0,w_1,\ldots,w_d$ such that $N(z)$ is an interpolating matrix polynomial:
	\[
	w_k = \frac{1}{\prod_{j\neq k} (z_j - z_k)};
	\]
	see, e.g., \cite[eq.~(3.2)]{BT04} and \cite[eq.~(3.4)]{VB15}. With this choice of the weights, we can immediately apply \cite[Thm.~4.8]{VB15}, which states that the matrix pencil $L(z) =  L_0 - z L_1$, where
	\[
	L_0 = \begin{bmatrix}
	z_1 C_0 & z_2 C_1 & \cdots & z_{d-1} C_{d-2} & z_\ell C_{d-1} + z_{d-1} \theta_d^{-1} C_d \\
	z_0 I & -z_2 \theta_1 I \\
	& \ddots & \ddots \\
	&        & z_{d-3} I & -z_{d-1}\theta_{d-2} I \\
	&   &                   & z_{d-2} I & -z_{d}\theta_{d-1} I 
	\end{bmatrix}
	\]
	and
	\[
	L_1 = \begin{bmatrix}
	C_0 &  C_1 & \cdots &  C_{d-2} & C_{d-1} + \theta_d^{-1} C_d \\
	I & - \theta_1 I \\
	& \ddots & \ddots  \\
	&        &  I & -\theta_{d-2} I \\
	&   &                   & I & -\theta_{d-1} I 
	\end{bmatrix}
	\]
	with $\theta_j = w_{j-1}/w_j$ for $j=1,\ldots,d$ is a strong linearization of $N(z)$. In other words, the nonlinear eigenvalues of $N(z)$ can be computed as the generalized eigenvalues of the matrix pair $(L_0,L_1)$. As long as $R_d$ does not also have a pole at a nonlinear eigenvalue $\lambda$ of $N(z)$, $\lambda$ will be a nonlinear eigenvalue of $R_d(z)$. 
	
	Finally, it is perhaps interesting to note that we have freedom to choose the nonzero weights $w_j$ and this choice will likely influence the numerical stability of the linearization $(L_0,L_1)$. This might be explored in some future work.

\end{document}